\documentclass{icmart}

\usepackage{amsmath,amssymb,mathrsfs}
\usepackage{amsrefs}
\usepackage[all]{xy}
\newdir{ >}{{}*!/-5pt/@{>}} 
\SelectTips{cm}{} 

\contact[rognes@math.uio.no]{Department of Mathematics,
	University of Oslo, Norway}

\newtheorem{theorem}{Theorem}[section]

\newtheorem{conjecture}[theorem]{Conjecture}
\theoremstyle{definition}
\newtheorem{definition}[theorem]{Definition}
\theoremstyle{remark}
\newtheorem{example}[theorem]{Example}
\newtheorem{remark}[theorem]{Remark}

\newcommand{\bC}{\mathbb{C}}
\newcommand{\bF}{\mathbb{F}}
\newcommand{\bG}{\mathbb{G}}

\newcommand{\bQ}{\mathbb{Q}}
\newcommand{\bW}{\mathbb{W}}
\newcommand{\bZ}{\mathbb{Z}}
\newcommand{\cC}{\mathcal{C}}
\newcommand{\cH}{\mathcal{H}}
\newcommand{\cJ}{\mathcal{J}}
\newcommand{\cS}{\mathcal{S}}
\newcommand{\sA}{\mathscr{A}}
\newcommand{\sO}{\mathscr{O}}
\newcommand{\sP}{\mathscr{P}}
\newcommand{\bn}{\mathbf{n}}

\newcommand{\Aut}{\operatorname{Aut}}
\newcommand{\Wh}{\operatorname{Wh}}
\newcommand{\PL}{\operatorname{PL}}
\newcommand{\Cat}{\operatorname{Cat}}
\newcommand{\Ext}{\operatorname{Ext}}
\newcommand{\Top}{\operatorname{Top}}
\newcommand{\Diff}{\operatorname{Diff}}
\newcommand{\Sp}{\operatorname{Sp}}
\newcommand{\Spec}{\operatorname{Spec}}

\newcommand{\trc}{\operatorname{trc}}
\newcommand{\colim}{\operatorname*{colim}}
\newcommand{\hocolim}{\operatorname*{hocolim}}
\newcommand{\cy}{\operatorname{cy}}
\newcommand{\rep}{\operatorname{rep}}
\newcommand{\longto}{\longrightarrow}
\newcommand{\CSp}{\cC\!\Sp^\Sigma}
\newcommand{\CSJ}{\cC\cS^{\cJ}}
\newcommand{\SJ}{S^{\cJ}}

\newcommand{\et}{\operatorname{\text{\'et}}}
\newcommand{\ff}{\operatorname{\text{\it ff}}}
\newcommand{\gp}{\operatorname{gp}}
\newcommand{\mot}{\operatorname{mot}}
\newcommand{\nr}{\operatorname{nr}}

\title{Algebraic $K$-theory of strict ring spectra}
\author{John Rognes}

\begin{document}
\begin{abstract}
We view strict ring spectra as generalized rings.  The study of their
algebraic $K$-theory is motivated by its applications to the automorphism
groups of compact manifolds.  Partial calculations of algebraic $K$-theory
for the sphere spectrum are available at regular primes, but we seek
more conceptual answers in terms of localization and descent properties.
Calculations for ring spectra related to topological $K$-theory suggest
the existence of a motivic cohomology theory for strictly commutative ring
spectra, and we present evidence for arithmetic duality in this theory.
To tie motivic cohomology to Galois cohomology we wish to spectrally
realize ramified extensions, which is only possible after mild forms of
localization. One such mild localization is provided by the theory of
logarithmic ring spectra, and we outline recent developments in this area.
\end{abstract}



\maketitle

\section{Strict ring spectra}

First, let $R$ be an abelian group.  Ordinary singular cohomology with
coefficients in $R$ is a contravariant homotopy functor that associates
to each based space $X$ a graded cohomology group $\tilde H^*(X; R)$.
It is stable, in the sense that there is a natural isomorphism $\tilde
H^*(X; R) \cong \tilde H^{*+1}(\Sigma X; R)$, and this implies that it
extends from the category of based spaces to the category of spectra.
The latter is a category of space-like objects, where the suspension
is invertible up to homotopy equivalence, and which has all colimits
and limits.  The extended cohomology functor becomes representable,
meaning that there is a spectrum $HR$, called the Eilenberg--Mac\,Lane
spectrum of $R$, and a natural isomorphism $\tilde H^*(X; R) \cong [X,
HR]_{-*}$, where $[X, HR]_{-n}$ is the group of homotopy classes of
morphisms $X \to \Sigma^n HR$.

Next, let $R$ be a ring.  Then the cohomology theory is multiplicative,
meaning that there is a bilinear cup product $\tilde H^*(X; R) \times
\tilde H^*(X; R) \to \tilde H^*(X; R)$.  This is also representable in the
category of spectra, by a morphism $\mu \colon HR \wedge HR \to HR$, where
$\wedge$ denotes the smash product of spectra.  With the modern models
for the category of spectra \cite{EKMM97}, \cite{HSS00}, \cite{MMSS01}
we may arrange that $\mu$ is strictly unital and associative, so that
$HR$ is a \emph{strict ring spectrum}.  Equivalent terms are $A_\infty$
ring spectrum, $S$-algebra, symmetric ring spectrum and orthogonal
ring spectrum.

If $R$ is commutative, then the cup product is graded commutative,
which at the representing level means that $\mu \tau \simeq \mu$,
where $\tau \colon HR \wedge HR \to HR \wedge HR$ denotes the twist
isomorphism.  In fact, we may arrange that $\mu$ is strictly commutative,
in the sense that $\mu\tau = \mu$ as morphisms of spectra, so that $HR$
is a \emph{strictly commutative ring spectrum}.  Equivalent phrases are
$E_\infty$ ring spectrum, commutative $S$-algebra, commutative symmetric
ring spectrum and commutative orthogonal ring spectrum.  This leads to
a compatible sequence of $\Sigma_k$-equivariant morphisms $E\Sigma_{k+}
\wedge HR^{\wedge k} \to HR$ for $k\ge0$.  At the represented level
these morphisms give rise to power operations in cohomology, including
Steenrod's operations $Sq^i$ for $R = \bF_2$ and $\beta^\epsilon P^i$
for $R = \bF_p$.

One now realizes that the Eilenberg--Mac\,Lane ring spectra $HR$
exist as special cases within a much wider class of ring spectra.
Each spectrum~$B$ represents a generalized cohomology theory $X \mapsto
\tilde B^*(X) = [X, B]_{-*}$ and a generalized homology theory $X \mapsto
\tilde B_*(X) = \pi_*(B \wedge X)$.  Examples of early interest include
the spectrum~$KU$ that represents complex topological $K$-theory, $KU^*(X)
= [X_+, KU]_{-*}$, and the spectrum~$MU$ that represents complex bordism,
$MU_*(X) = \pi_*(MU \wedge X_+)$.  A fundamental example is given by
the sphere spectrum~$S$, which is the image of the based space $S^0$
under the stabilization functor from spaces to spectra.  It represents
stable cohomotopy $\pi^*_S(X) = \tilde S^*(X)$ and stable homotopy
$\pi^S_*(X) = \tilde S_*(X)$.  Each of these three examples, $KU$, $MU$
and $S$, is naturally a strictly commutative ring spectrum, representing a
multiplicative cohomology theory with power operations, etc.  Furthermore,
there are interesting multiplicative morphisms connecting these ring
spectra to the Eilenberg--Mac\,Lane ring spectra previously considered,
as in the diagram
$$
\xymatrix{
& & KU \ar[r] & H\bQ \\
S \ar[r] & MU \ar[r] & ku \ar[r] \ar[u] & H\bZ \ar[u] \rlap{\,,}
}
$$
where $ku = KU[0,\infty)$ denotes the connective cover of $KU$.

By placing the class of traditional rings inside the wider realm
of all strict ring spectra, a new world of possibilities opens up.
Following Waldhausen \cite{MP89}*{p.~xiii} we may refer to strict ring
spectra as ``brave new rings''.  If we think in algebro-geometric terms,
where commutative rings appear as the rings of functions on pieces
of geometric objects, then strictly commutative ring spectra are the
functions on affine pieces of brave new geometries, more general than
those realized by ordinary schemes.

How vast is this generalization?  In the case of connective ring
spectra~$B$, i.e., those with $\pi_i(B) = 0$ for $i$ negative, there
is a natural ring spectrum morphism $B \to H\pi_0(B)$ that induces
an isomorphism on~$\pi_0$.  This behaves for many purposes like a
topologically nilpotent extension, and in geometric terms, $B$ can be
viewed as the ring spectrum of functions on an infinitesimal thickening
of $\Spec \pi_0(B)$.

This infinitesimal thickening can be quite effectively controlled
in terms of diagrams of Eilenberg--Mac\,Lane spectra associated with
simplicial rings.  The Hurewicz map $B \cong S \wedge B \to H\bZ \wedge
B$ is $1$-connected, and there is an equivalence $H\bZ \wedge B \simeq
HR_\bullet$ for some simplicial ring~$R_\bullet$.  The square
$$
\xymatrix{
S \wedge S \wedge B \ar[r] \ar[d] & H\bZ \wedge S \wedge B \ar[d] \\
S \wedge H\bZ \wedge B \ar[r] & H\bZ \wedge H\bZ \wedge B
}
$$
induces a $2$-connected map from $B \cong S \wedge S \wedge B$ to
the homotopy pullback.  More generally, for each $n\ge1$ there is
an $n$-dimensional cubical diagram that induces an $n$-connected map
from the initial vertex $B$ to the homotopy limit of the remainder of
the cube, and the terms in that remainder have the form $HR_\bullet$
for varying simplicial rings~$R_\bullet$.  Dundas \cite{Dun97} used a
clever strengthening of this statement to prove that relative algebraic
$K$-theory is $p$-adically equivalent to relative topological cyclic
homology for morphisms $A \to B$ of connective strict ring spectra,
under the assumption that $\pi_0(A) \to \pi_0(B)$ is a surjection with
nilpotent kernel.  He achieved this by reducing to the analogous statement
for homomorphisms $R_\bullet \to T_\bullet$ of simplicial rings, which
had been established earlier by McCarthy \cite{McC97}.  This confirmed
a conjecture of Goodwillie \cite{Goo91}, motivated by a similar result
for rational $K$-theory and (negative) cyclic homology \cite{Goo86}.

How about the case of non-connective ring spectra?  Those that arise as
homotopy fixed points $B^{hG} = F(EG_+, B)^G$ for a group action may be
viewed as the functions on an orbit stack for the induced $G$-action on
the geometry associated to~$B$.  Those that arise as smashing Bousfield
localizations $L_E B$, with respect to a homology theory $E_*$, may
be viewed as open subspaces in a finer topology than the one derived
from the Zariski topology on $\Spec \pi_0(B)$ \cite{Rog08}*{\S9.3}.
In the general case the connection to classical geometry is less clear.

\section{Automorphisms of manifolds}

Why should we be interested in brave new rings and their ring-theoretic
invariants, like algebraic $K$-theory \cite{Qui75}, \cite{EKMM97}*{Ch.~VI},
other than for the sake of generalization?  One good justification
comes from the tight connection between the geometric topology of
high-dimensional manifolds and the algebraic $K$-theory of strict ring
spectra.  This connection is given by the higher simple-homotopy theory
initiated by Hatcher \cite{Hat75}, which was fully developed by Waldhausen
in the context of his algebraic $K$-theory of spaces \cite{Wal85} and
the stable parametrized $h$-cobordism theorem \cite{WJR13}.  On the
geometric side, this theory concerns the fundamental problem of finding
a parametrized classification of high-dimensional compact manifolds,
up to homeomorphism, piecewise-linear homeomorphism or diffeomorphism,
as appropriate for the respective geometric category.  The set of
path components of the resulting moduli space corresponds to the set
of isomorphism classes of such manifolds, and each individual path
component is a classifying space for the automorphism group $\Aut(X)$
of a manifold~$X$ in the respective isomorphism class.

This parametrized classification is finer than the one provided by
the Browder--Novikov--Sullivan--Wall surgery theory \cite{Bro72},
\cite{Wll70}, which classifies manifolds up to $h$-cobordism (or
$s$-cobordism), and whose associated moduli space has path components
that classify the block automorphism groups of manifolds, rather than
their actual automorphism groups.  The difference between these two
classifications is controlled by the space $H(X)$ of $h$-cobordisms $(W;
X, Y)$ with a given manifold $X$ at one end.  Here $W$ is a compact
manifold with $\partial W = X \cup Y$, and the inclusions $X \to W$
and $Y \to W$ are homotopy equivalences.  More precisely, there is
one $h$-cobordism space $H^{\Cat}(X)$ for each of the three flavors of
manifolds mentioned, namely $\Cat = \Top$, $\PL$ or~$\Diff$.

The original $h$-cobordism theorem enumerates the isomorphism classes of
$h$-cobordisms with $X$ at one end, i.e., the set $\pi_0 H(X)$ of path
components of the $h$-cobordism space of $X$, in terms of an algebraic
$K$-group of the integral group ring $\bZ[\pi]$, where $\pi = \pi_1(X)$
is the fundamental group of $X$.  One defines the Whitehead group as the
quotient $\Wh_1(\pi) = K_1(\bZ[\pi])/(\pm\pi)$, and associates a Whitehead
torsion class $\tau(W; X, Y) \in \Wh_1(\pi)$ to each $h$-cobordism on $X$.

\begin{theorem}[Smale \cite{Sma62}, Barden, Mazur, Stallings]
Let $X$ be a compact, connected $n$-manifold with $n\ge5$ and $\pi =
\pi_1(X)$.  The Whitehead torsion defines a bijection
$$
\pi_0 H(X) \cong \Wh_1(\pi) \,.
$$
\end{theorem}

These constructions involve using Morse functions or triangulations
to choose a relative CW complex structure on the pair $(W, X)$, or
equivalently, a $\pi$-equivariant relative CW complex structure on
the pair of universal covers $(\tilde W, \tilde X)$, and to study the
associated cellular chain complex $C_*(\tilde W, \tilde X)$ of free
$\bZ[\pi]$-modules.  This works fine as long as one is only concerned
with a classification of $h$-cobordisms up to isomorphism, but for the
parametrized problem, i.e., the study of the full homotopy type of $H(X)$,
the passage from a CW complex structure to the associated cellular chain
complex loses too much information.  One should remember the actual
attaching maps from the boundaries of cells to the preceding skeleta, not
just their degrees. In a stable range this amounts to working with maps
from spheres to spheres and coefficients in the sphere spectrum~$S$,
rather than with degrees and coefficients in the integers~$\bZ$.
Likewise, the passage to the $\pi$-equivariant universal cover $\tilde
X \to X$ should be replaced to a passage to a $G$-equivariant principal
fibration $P \to X$, where $P$ is contractible and $G \simeq \Omega X$
is a topological group that is homotopy equivalent to the loop space
of $X \simeq BG$.  To sum up, the parametrized analog of the Whitehead
torsion must take values in a Whitehead space that is built from the
algebraic $K$-theory of the spherical group ring $S[G] = S \wedge G_+$,
a strict ring spectrum, rather than that of its discrete reduction,
the integral group ring $\pi_0(S[G]) \cong \bZ[\pi]$.

Waldhausen's algebraic $K$-theory of spaces, traditionally denoted
$A(X)$, was first introduced without reference to strict ring spectra
\cite{Wal78}, but can be rewritten as the algebraic $K$-theory $A(X)
= K(S[G])$ of the strict ring spectrum $S[G]$, cf.~\cite{Wal84} and
\cite{EKMM97}*{Ch.~VI}.  This point of view is convenient for the
comparison of algebraic $K$-theory with other ring-theoretic invariants.

In the case of differentiable manifolds, the Whitehead space
$\Wh^{\Diff}(X)$ is defined to sit in a split homotopy fiber sequence
of infinite loop spaces
$$
\Omega^\infty(S \wedge X_+)
	\overset{\iota}\longto K(S[G]) \longto \Wh^{\Diff}(X) \,.
$$
In the topological case there is a homotopy fiber sequence
of infinite loop spaces
$$
\Omega^\infty(K(S) \wedge X_+)
	\overset{\alpha}\longto K(S[G]) \longto \Wh^{\Top}(X) \,,
$$
where $\alpha$ is known as the assembly map.  The piecewise-linear
Whitehead space is the same as the topological one.
The stable parametrized $h$-cobordism theorem reads as follows.

\begin{theorem}[Waldhausen--Jahren--Rognes \cite{WJR13}*{Thm.~0.1}]
Let $X$ be a compact $\Cat$ manifold, for $\Cat = \Top$, $\PL$ or $\Diff$.
There is a natural homotopy equivalence
$$
\cH^{\Cat}(X) \simeq \Omega \Wh^{\Cat}(X) \,,
$$
where $\cH^{\Cat}(X) = \colim_k H^{\Cat}(X \times I^k)$ is the stable
$\Cat$ $h$-cobordism space of $X$.
\end{theorem}

When combined with connectivity results about the dimensional
stabilization map $H(X) = H^{\Cat}(X) \to \cH^{\Cat}(X)$, and here the
main result is Igusa's stability theorem for smooth pseudoisotopies
\cite{Igu88}, knowledge of $K(S)$ and $K(S[G])$ gives good general
results on the $h$-cobordism space $H(X)$ and the automorphism
group $\Aut(X)$ of a high-dimensional manifold~$X$.

\begin{example}
When $G$ is trivial, so that $S[G] = S$ and $X$ is contractible, the
$\pi_0$-isomorphism and rational equivalence $S \to H\bZ$ induces a
rational equivalence $K(S) \to K(\bZ)$.  Here $\pi_* K(\bZ) \otimes \bQ$
was computed by Borel \cite{Bor74}, so
$$
\pi_i \Wh^{\Diff}(*) \otimes \bQ \cong
\begin{cases}
\bQ & \text{for $i = 4k+1 \ne 1$,} \\
0 & \text{otherwise.}
\end{cases}
$$
For $X = D^n$, Farrell--Hsiang \cite{FH78} used this to
show that
$$
\pi_i \Diff(D^n) \otimes \bQ \cong
\begin{cases}
\bQ & \text{for $i = 4k-1$, $n$ odd,} \\
0 & \text{otherwise,}
\end{cases}
$$
for $i$ up to approximately $n/3$, where $\Diff(D^n)$ denotes the group
of self-diffeo\-morphisms of $D^n$ that fix the boundary.  For instance,
$\pi_3 \Diff(D^{13})$ is rationally nontrivial.  By contrast, the group
$\Top(D^n)$ of self-homeomorphisms of $D^n$ that fix the boundary
is contractible.  Similar results follow for $n$-manifolds that are
roughly $n/3$-connected.
\end{example}

The case of spherical space forms, when $G$ is finite with periodic
cohomology, has been studied by Hsiang--Jahren \cite{HJ82}.  For closed,
non-positively curved manifolds~$X$, Farrell--Jones \cite{FJ91} showed
that $\Wh^{\Diff}(X)$ can be assembled from copies of $\Wh^{\Diff}(*)$
and $\Wh^{\Diff}(S^1)$, indexed by the points and the closed geodesics
in~$X$, respectively.  These correspond to the special cases $G$ trivial
and $G$ infinite cyclic, respectively, so $K(S)$ and $K(S[\bZ])$ are
of fundamental importance for the parametrized classification of this
large class of Riemannian manifolds.  In this paper we shall focus on the
case of $K(S)$, but see Hesselholt's paper \cite{Hes09} for the case of
$K(S[\bZ])$, and see Weiss--Williams \cite{WW01} for a detailed survey
about automorphisms of manifolds and algebraic $K$-theory.

\begin{remark}
More recent papers of Madsen--Weiss \cite{MW07}, Berglund--Madsen
\cite{BeM:14} and Galatius--Randal-Williams \cite{GRW:12} give precise
results about automorphism groups of manifolds of a fixed even dimension
$n = 2d \ne 4$, at the expense of first forming a connected sum with
many copies of $S^d \times S^d$.  The latter results are apparently not
closely related to the algebraic $K$-theory of strict ring spectra.
\end{remark}

\section{Algebraic $K$-theory of the sphere spectrum}

We can strengthen the rational results about $A(X) = K(S[G])$,
$\Wh^{\Diff}(X)$ and $\Diff(X)$ to integral results, or more precisely,
to $p$-adic integral results for each prime $p$.  From here on it will
be convenient to think of algebraic $K$-theory as a spectrum-valued
functor, and likewise for the Whitehead theories, so that there are
homotopy cofiber sequences of spectra
\begin{align*}
S \wedge X_+ &\overset{\iota}\longto K(S[G]) \longto \Wh^{\Diff}(X) \\
K(S) \wedge X_+ &\overset{\alpha}\longto K(S[G])
	\longto \Wh^{\Top}(X) \,,
\end{align*}
and the first one is naturally split.

A key tool for this study is the cyclotomic trace map
$\trc \colon K(B) \to TC(B; p)$
from algebraic $K$-theory to the topological cyclic homology of
B{\"o}kstedt--Hsiang--Madsen \cite{BHM93}.  The latter invariant of the
strict ring spectrum $B$ can sometimes be calculated by analyzing the
$S^1$-equivariant homotopy type of the topological Hochschild homology
spectrum $THH(B)$.  Its power is illustrated by the following previously
mentioned theorem.

\begin{theorem}[Dundas \cite{Dun97}]
Let $B$ be a connective strict ring spectrum.  The square
$$
\xymatrix{
K(B) \ar[r] \ar[d]_-{\trc} & K(\pi_0(B)) \ar[d]^-{\trc} \\
TC(B; p) \ar[r] & TC(\pi_0(B); p)
}
$$
becomes homotopy Cartesian upon $p$-completion.
\end{theorem}

In the basic case $B = S$, when $K(S) \simeq S \vee \Wh^{\Diff}(*)$
determines $\Diff(D^n)$ for large~$n$, this square takes the form below.
Three of the four corners are quite well understood, but for widely
different reasons.
$$
\xymatrix{
K(S) \ar[r] \ar[d]_-{\trc} & K(\bZ) \ar[d]^-{\trc} \\
TC(S; p) \ar[r] & TC(\bZ; p) \rlap{\,.}
}
$$
These reasons were tied together by the author for $p=2$ in \cite{Rog02},
and for $p$ an odd regular prime in \cite{Rog03}, to compute the mod~$p$
cohomology
$$
H^*(K(S); \bF_p) \cong \bF_p \oplus H^*(\Wh^{\Diff}(*); \bF_p)
$$
as a module over the Steenrod algebra $\sA$ of stable mod~$p$ cohomology
operations.  This sufficed to determine the $E_2$-term of the Adams
spectral sequence
$$
E_2^{s,t} = \Ext_{\sA}^{s,t}(H^*(K(S); \bF_p), \bF_p)
	\Longrightarrow \pi_{t-s} K(S)_p
$$
in a large range of degrees, and to determine the homotopy groups $\pi_i
K(S)_p$ and $\pi_i \Wh^{\Diff}(*)_p$ in a smaller range of degrees.

The structure of the algebraic $K$-theory of the integers, $K(\bZ)$,
was predicted by the Lichtenbaum--Quillen conjectures \cite{Qui75}, which
were confirmed for $p=2$ by Voevodsky \cite{Voe03} with contributions by
Rognes--Weibel \cite{RW00}, and for $p$ odd by Voevodsky \cite{Voe11}
with contributions by Rost and Weibel.  For $p=2$ or $p$ a regular odd
prime, this led to a $p$-complete description of the spectrum $K(\bZ)$
in terms of topological $K$-theory spectra, which in turn led to an
explicit description of the spectrum cohomology $H^*(K(\bZ); \bF_p)$
as an $\sA$-module.

The topological cyclic homology of the integers, $TC(\bZ; p)$, was
computed for odd primes $p$ by B{\"o}kstedt--Madsen \cite{BM94},
\cite{BM95} and for $p=2$ by the author \cite{Rog98}, \cite{Rog99a},
\cite{Rog99b}, \cite{Rog99c}, in papers that start with
knowledge of the mod~$p$ homotopy of the $S^1$-spectrum $THH(\bZ)$ and
inductively determine the mod~$p$ homotopy of the $C_{p^n}$-fixed points
$THH(\bZ)^{C_{p^n}}$ for $n\ge1$.  It is then possible to recognize
the $p$-completed spectrum level structure by comparisons with known
models, using \cite{Rog93} for $p$ odd, and to obtain the $\sA$-module
$H^*(TC(\bZ; p); \bF_p)$ from this.

The topological cyclic homology of the sphere spectrum, $TC(S; p)$, was
determined in the original paper \cite{BHM93}.  There is an equivalence
of spectra $TC(S; p) \simeq S \vee \Sigma \bC P^\infty_{-1}$ after
$p$-completion, where $\Sigma \bC P^\infty_{-1}$ is the homotopy fiber of
the dimension-shifting $S^1$-transfer map $t \colon \Sigma \bC P^\infty_+
\to S$.  The mod~$p$ cohomology $H^*(\Sigma \bC P^\infty_{-1}; \bF_p)$
is well known as an $\sA$-module.  For $p=2$ it is cyclic with
$$
H^*(\Sigma \bC P^\infty_{-1}; \bF_2) \cong \Sigma^{-1} \sA/C \,,
$$
where the ideal $C \subset \sA$ is generated by the admissible $Sq^I$
where $I = (i_1, \dots, i_n)$ with $n\ge2$ or $I = (i)$ with $i$ odd.
The determination of the homotopy groups of~$TC(S; p)$ is of comparable
difficulty to the computation of the homotopy groups of~$S$, due to our
extensive knowledge about the attaching maps in the usual CW spectrum
structure on $\Sigma \bC P^\infty_{-1}$, cf.~Mosher \cite{Mos68}.

The linearization map $TC(S; p) \to TC(\bZ; p)$ is only partially
understood \cite{KR97}, but for $p$ regular the cyclotomic trace
map $K(\bZ) \to TC(\bZ; p)$ can be controlled by an appeal
to global Tate--Poitou duality \cite{Tat63}*{Thm.~3.1}, see
\cite{Rog03}*{Prop.~3.1}.  This leads to the following conclusion for
$p=2$.  See \cite{Rog03}*{Thm.~5.4} for the result at odd regular primes.

\begin{theorem}[\cite{Rog02}*{Thm.~4.5}]
The mod~$2$ cohomology of the spectrum $\Wh^{\Diff}(*)$ is given
by the unique non-trivial extension of $\sA$-modules
$$
0 \to \Sigma^{-2} C/\sA(Sq^1, Sq^3)
	\longto H^*(\Wh^{\Diff}(*); \bF_2)
	\longto \Sigma^3 \sA/\sA(Sq^1, Sq^2) \to 0 \,.
$$
\end{theorem}

Using the Adams spectral sequence and related methods, the author
obtained the following explicit calculations.  Less complete information,
in a larger range of degrees, is provided in the cited references.
Previously, B{\"o}kstedt--Waldhausen \cite{BW87}*{Thm.~1.3} had computed
$\pi_i \Wh^{\Diff}(*)$ for $i\le3$.

\begin{theorem}[\cite{Rog02}*{Thm.~5.8}, \cite{Rog03}*{Thm.~4.7}]
The homotopy groups of $\Wh^{\Diff}(*)$ in degrees $i\le18$ are as
follows, modulo $p$-power torsion for irregular primes~$p$.
$$
\begin{tabular}{ c | c c c c c c c c c c }
$i$ & $0$ & $1$ & $2$ & $3$ & $4$
	& $5$ & $6$ & $7$ & $8$
	& $9$ \\
\hline
$\pi_i \Wh^{\strut\Diff}(*)$ & $0$ & $0$ & $0$ & $\bZ/2$ & $0$ 
	& $\bZ$ & $0$ & $\bZ/2$ & $0$
	& $\bZ \oplus \bZ/2$
\end{tabular}
$$
$$
\begin{tabular}{ c | c c c c c }
$i$ & $10$ & $11$ & $12$
	& $13$ & $14$ \\
\hline
$\pi_i \Wh^{\strut\Diff}(*)$ & $\bZ/8 \oplus (\bZ/2)^2$ & $\bZ/6$ & $\bZ/4$
	& $\bZ$ & $\bZ/36 \oplus \bZ/3$
\end{tabular}
$$
$$
\begin{tabular}{ c | c c c c }
$i$ & $15$ & $16$
	& $17$ & $18$ \\
\hline
$\pi_i \Wh^{\strut\Diff}(*)$ & $(\bZ/2)^2$ & $\bZ/24 \oplus \bZ/2$
	& $\bZ \oplus (\bZ/2)^2$ & $\bZ/480 \oplus (\bZ/2)^3$
\end{tabular}
$$
\end{theorem}

\begin{example}
For $X = D^n$ with $n$ sufficiently large, it follows that $\pi_{4p-4}
\Diff(D^n)$ or $\pi_{4p-4} \Diff(D^{n+1})$ contains an element of
order~$p$, for each regular $p\ge5$, and that $\pi_9 \Diff(D^n)$
or $\pi_9 \Diff(D^{n+1})$ contains an element of order~$3$, see
\cite{Rog03}*{Thm.~6.4}.  To get more precise results one needs to
investigate the canonical involution on $\Wh^{\Diff}(*)$ and apply
Weiss--Williams \cite{WW88}*{Thm.~A}.
\end{example}

\begin{remark}
It would be interesting to extend these results to irregular primes.
Dwyer--Mitchell \cite{DM98} described the spectrum $K(\bZ)_p$ in terms of
the $p$-primary Iwasawa module of the rationals.  It should be possible to
turn this into a description of the $\sA$-module $H^*(K(\bZ); \bF_p)$.
Next one must control the cyclotomic trace map $K(\bZ) \to TC(\bZ;
p)$, or the closely related completion map $K(\bZ) \to K(\bZ_p)$,
whose behavior is governed by special values of $p$-adic $L$-functions,
cf.~Soul{\'e} \cite{Sou81}*{Thm.~3}.
\end{remark}

\section{Algebraic $K$-theory of topological $K$-theory}

The calculations reviewed in the previous section extracted detailed
information about $\pi_* K(S) \cong \pi_*(S) \oplus \pi_* \Wh^{\Diff}(*)$
from our knowledge of $\pi_*(\bC P^\infty_{-1})$.  However, this
understanding was not presented to us in as conceptual a way as the
understanding we have of $K(\bZ)$, say in terms of Quillen's localization
sequence
$$
K(\bF_p) \longto K(\bZ) \longto K(\bZ[1/p])
$$
and the {\'e}tale descent property
$$
\pi_i K(\bZ[1/p])_p \overset{\cong}\longto
	K^{\et}_i(\bZ[1/p]; \bZ_p)
$$
for $i>0$, cf.~\cite{Qui73}*{\S5} and~\cite{Qui75}*{\S9}.  It would
be desirable to have a similarly conceptual understanding of $K(S)_p$
in terms of a comparison with $K(B)_p$ for suitably local strict ring
spectra~$B$, a descent property describing $K(B)_p$ as a homotopy limit
of $K(C)_p$ for appropriate extensions $B \to C$, and a simple description
of $K(\Omega)_p$ for a sufficiently large such extension $B \to \Omega$.

To explore this problem, we first simplify the number theory involved
by working with the $p$-adic integers $\bZ_p$ in place of the rational
integers $\bZ$, and then seek a conceptual understanding of $K(B)_p$
for some of the strictly commutative ring spectra~$B$ that are closest
to $H\bZ_p$, namely the $p$-complete connective complex $K$-theory
spectrum $ku_p$ and its Adams summand $\ell_p$.  Here $\pi_*(ku_p) =
\bZ_p[u]$ and $\pi_*(\ell_p) = \bZ_p[v_1]$, with $|u| = 2$ and $|v_1|
= 2p-2$.  Let $KU_p$ and $L_p$ denote the associated periodic spectra,
with $\pi_*(KU_p) = \bZ_p[u^{\pm1}]$ and $\pi_*(L_p) = \bZ_p[v_1^{\pm1}]$.
There are multiplicative morphisms
\begin{equation} \label{eq:SphiHZp}
\xymatrix{
& L_p \ar[r]^-{\phi} & KU_p \\
S_p \ar[r] & \ell_p \ar[r]^-{\phi} \ar[u]
	& ku_p \ar[r] \ar[u] & H\bZ_p
}
\end{equation}
of strictly commutative ring spectra, where $\phi_*(v_1) = u^{p-1}$.
The group $\Delta \cong \bF_p^\times$ of $p$-adic roots of unity
acts by Adams operations on $KU_p$, and $\phi \colon L_p \to KU_p$ is a
$\Delta$-Galois extension in the sense of \cite{Rog08}*{p.~3}.

\begin{definition}
Let $V(1) = S \cup_p e^1 \cup_{\alpha_1} e^{2p-1} \cup_p e^{2p}$ be the
type~$2$ Smith--Toda complex, defined as the mapping cone of the Adams
self-map $v_1 \colon \Sigma^{2p-2} S/p \to S/p$ of the mod~$p$ Moore spectrum
$S/p = S \cup_p e^1$.  It is a ring spectrum up to homotopy for $p\ge5$,
which we now assume.  We write $V(1)_* B = \pi_*(V(1) \wedge B)$
for the ``mod~$p$ and~$v_1$ homotopy'' of any spectrum $B$.  It is
naturally a module over the polynomial ring $\bF_p[v_2]$,
where $v_2 \in \pi_{2p^2-2} V(1)$.
\end{definition}

The mod~$p$ and~$v_1$ homotopy of the topological cyclic homology of the
connective Adams summand $\ell$ was computed by Ausoni and the author,
by starting with knowledge of $V(1)_* THH(\ell)$ from \cite{MS93}
and inductively determining the mod~$p$ and~$v_1$ homotopy of the
fixed points $THH(\ell)^{C_{p^n}}$ for $n\ge1$.  The calculations were
later extended to the full connective complex $K$-theory spectrum $ku$
by Ausoni.  To avoid introducing too much notation, we only describe the
most striking features of the answers, referring to the original papers
for more precise statements.

\begin{theorem}[Ausoni--Rognes \cite{AR02}*{Thm.~0.3, Thm.~0.4}]
\label{thm:TCellp}
$V(1)_* TC(\ell; p)$ is a finitely generated free
$\bF_p[v_2]$-module on $4(p+1)$ generators, which are located in degrees
$-1 \le * \le 2p^2+2p-2$.
There is an exact sequence of $\bF_p[v_2]$-modules
$$
0 \to \Sigma^{2p-3} \bF_p \longto V(1)_* K(\ell_p) \overset{\trc}\longto
	V(1)_* TC(\ell; p) \longto \Sigma^{-1} \bF_p \to 0 \,.
$$
\end{theorem}

\begin{theorem}[Ausoni \cite{Aus10}*{Thm.~7.9, Thm.~8.1}]
$V(1)_* TC(ku; p)$ is a finitely generated free
$\bF_p[v_2]$-module on $4(p-1)(p+1)$ generators, which are located in degrees
$-1 \le * \le 2p^2+2p-2$.
There is an exact sequence of $\bF_p[v_2]$-modules
$$
0 \to \Sigma^{2p-3} \bF_p \longto V(1)_* K(ku_p) \overset{\trc}\longto
	V(1)_* TC(ku; p) \longto \Sigma^{-1} \bF_p \to 0 \,,
$$
and the natural map $K(\ell_p) \to K(ku_p)^{h\Delta}$ is a $p$-adic
equivalence.
\end{theorem}

Blumberg--Mandell \cite{BM08} constructed homotopy cofiber sequences
\begin{equation}
\label{eq:Klocell}
K(\bZ_p) \longto K(\ell_p) \longto K(L_p)
\end{equation}
and
$K(\bZ_p) \to K(ku_p) \to K(KU_p)$,
which lead to calculations of $V(1)_* K(L_p)$ and $V(1)_*
K(KU_p)$, cf.~\cite{Aus10}*{Thm.~8.3}.  The natural map $K(L_p) \to
K(KU_p)^{h\Delta}$ is also a $p$-adic equivalence, which confirms the
{\'e}tale descent property for algebraic $K$-theory in this particular
case.

\begin{remark}
The examples discussed above are the case $n=1$ of a series of
approximations to $S$ associated with the Lubin--Tate spectra $E_n$,
with coefficient rings $\pi_* E_n = \bW\bF_{p^n}[[u_1, \dots,
u_{n-1}]][u^{\pm1}]$, which are known to be strictly commutative ring
spectra by the Goerss--Hopkins--Miller obstruction theory \cite{GH04}.
There are multiplicative morphisms
$$
\xymatrix{
\hat L_n S \ar[r] & \hat L_n E(n) \ar[r]^-{\phi} & E_n \\
L_n S \ar[u] \ar[r] & BP\langle n\rangle \ar[u] \ar[r] & e_n \ar[u] \ar[r]
	& H\pi_0(e_n) \,,
}
$$
where $L_n$ and $\hat L_n$ denote Bousfield localization with respect
to the Johnson--Wilson spectrum~$E(n)$ and the Morava $K$-theory
spectrum~$K(n)$, respectively, and $BP\langle n\rangle$ is the truncated
Brown--Peterson spectrum.  The $n$-th extended Morava stabilizer group
$\bG_n$ acts on $E_n$, and $\hat L_n E(n) \to E_n$ is an $H$-Galois
extension for $H \cong \bF_{p^n}^\times \rtimes \bZ/n$.  We write $e_n$
for the connective cover of~$E_n$.

The algebraic $K$-theory computations above provide evidence for the
chromatic redshift conjecture, see \cite{AR02}*{p.~7} and~\cite{AR08},
predicting that the algebraic $K$-theory $K(B)$ of a purely $v_n$-periodic
strictly commutative ring spectrum~$B$, such as $E_n$, is purely
$v_{n+1}$-periodic in sufficiently high degrees.
\end{remark}

\section{Motivic truncation and arithmetic duality}

The proven Lichtenbaum--Quillen conjectures subsume a spectral sequence
$$
E^2_{s,t} = H_{\et}^{-s}(R; \bZ_p(t/2))
	\Longrightarrow \pi_{s+t} K(R)_p \,,
$$
which converges for reasonable $R$ and $s+t$ sufficiently large.
Here $H_{\et}^*$ denotes {\'e}tale cohomology, $R$ is a commutative
$\bZ[1/p]$-algebra, and $\bZ_p(t/2) = \pi_t(KU_p)$ is $\bZ_p(m)$
when $t = 2m$ is even, and $0$ otherwise.  For instance, we may take $R =
\sO_F[1/p]$ to be the ring of $p$-integers in a number field~$F$, or $R$
may be a $p$-adic field, i.e., a finite extension of $\bQ_p$.

The proven Beilinson--Lichtenbaum conjectures, cf.~\cite{SV00} and
\cite{GL01}, provide a more precise convergence statement.  For each
field $F$ containing $1/p$ there is a spectral sequence
$$
E^2_{s,t} = H_{\mot}^{-s}(F; \bZ(t/2))
	\Longrightarrow \pi_{s+t} K(F) \,,
$$
converging in all degrees, and similarly with mod~$p$ coefficients.
Here $H_{\mot}^*$ denotes motivic cohomology, which satisfies
\begin{equation} \label{eq:BLmotet}
H_{\mot}^r(F; \bZ/p(m)) \cong
\begin{cases}
H_{\et}^r(F; \bZ/p(m)) & \text{for $0 \le r \le m$,} \\
0 & \text{otherwise.}
\end{cases}
\end{equation}
In terms of Bloch's higher Chow groups \cite{Blo86}, the vanishing
of these groups for $r>m$ expresses the fact that there are no
codimension~$r$ subvarieties of affine $m$-space over $\Spec F$.
Conversely,
\begin{equation} \label{eq:BLetmot}
H_{\et}^r(F; \bZ/p(*)) \cong
	v_1^{-1} H_{\mot}^r(F; \bZ/p(*))
\end{equation}
with $v_1 \in H_{\mot}^0(F; \bZ/p(p-1))$.  We refer to
the aspects~\eqref{eq:BLmotet} and~\eqref{eq:BLetmot} of the
Beilinson--Lichtenbaum conjectures as the \emph{motivic truncation
property} for the field~$F$.

The following prediction expresses a similar conceptual description
of $K(B)$ for some strictly commutative ring spectra, and should in
particular apply for $B = \ell_p$, $L_p$, $ku_p$ and~$KU_p$.

\begin{conjecture} \label{conj:hmotk}
For purely $v_1$-periodic strictly commutative ring spectra $B$
there is a spectral sequence
$$
E^2_{s,t} = H_{\mot}^{-s}(B; \bF_{p^2}(t/2))
	\Longrightarrow V(1)_{s+t} K(B) \,,
$$
converging for $s+t$ sufficiently large.
\end{conjecture}

Here $H_{\mot}^*$ denotes a currently undefined form of motivic
cohomology for strictly commutative ring spectra.  The coefficient
$\bF_{p^2}(t/2)$ may be interpreted as $V(1)_t E_2$, where $E_2$ is
the Lubin--Tate ring spectrum \cite{GH04} with $\pi_* E_2
= \bW\bF_{p^2}[[u_1]][u^{\pm1}]$.

More generally one might consider purely $v_n$-periodic ring spectra $B$,
replace $V(1)$ by any type~$n+1$ finite spectrum $F(n+1)$, see \cite{HS98},
and replace $V(1)_t E_2$ and $V(1)_{s+t} K(B)$ by $F(n+1)_t E_{n+1}$
and $F(n+1)_{s+t} K(B)$, respectively.

\begin{example}
Based on the detailed calculations behind Theorem~\ref{thm:TCellp}, it is
fairly evident that the $E^2$-term for the spectral sequence conjectured
to converge to $V(1)_* K(\ell_p)$ will be concentrated in the four
columns $-3 \le s \le 0$, and that the free $\bF_p[v_2]$-module generators
are located in the groups $H_{\mot}^r(\ell_p; \bF_{p^2}(m))$ where
$0 \le r \le 3$ and $r \le m < r + p^2 + p - 1$.  This presumes that the
spectral sequence collapses at the $E^2$-term, for $p\ge5$.  In addition,
there is a sporadic copy of $\bF_p$ in $V(1)_{2p-3} K(\ell_p)$.

The class $v_2 \in V(1)_{2p^2-2} K(\ell_p)$ is represented in bidegree
$(s, t) = (0, 2p^2-2)$, corresponding to $(r, m) = (0, p^2-1)$.
The presence of $\bF_p[v_2]$-module generators in the range $r + p^2-1 \le m$
shows that $H_{\mot}^r(\ell_p; \bF_{p^2}(*))$ is \emph{not} isomorphic to
$v_2^{-1} H_{\mot}^r(\ell_p; \bF_{p^2}(*))$ in several bigradings $(r,
m)$ with $r \le m < r + p$.  In other words, the motivic truncation
property fails for $\ell_p$.  However, this is to be expected, since
$\ell_p$ has the residue ring spectrum $H\bZ_p$ and should not behave
as a field.
\end{example}

\begin{example}
Turning instead to $V(1)_* K(L_p)$, as determined from $V(1)_* K(\bZ_p)$
and $V(1)_* K(\ell_p)$ by the homotopy cofiber sequence~\eqref{eq:Klocell},
free $\bF_p[v_2]$-module generators for the $E^2$-term in
Conjecture~\ref{conj:hmotk} would be concentrated in the
groups $H^r_{\mot}(L_p; \bF_{p^2}(m))$ with $0 \le r \le 3$ and $r \le
m < r + p^2 - 1$.  The $(s, t)$-bidegrees of these $4p+4$ generators is
displayed for $p=5$ in Figure~\ref{fig:E2KLp}, lying in a fundamental
domain in the shape of a parallelogram, of width~$3$ and height $2p^2-2$.
In addition, there are sporadic copies of $\bF_p$ in $V(1)_{2p-3} K(L_p)$
and $V(1)_{2p-2} K(L_p)$.

\begin{figure}[t]
$$
\xymatrix@!C@!R@R-30pt{
\ar@{..}[3,3] \ar@{..}[24,0] \\
\\
& \bullet \bullet & & & 2p^2 \\
& & \bullet & \ar@{..}[24,0] \\
\\
\\
\\
& \bullet & \bullet & & - \\
\\
\\
\\
\\
& \bullet & \bullet & & - \\
\\
\\
\\
\\
& \bullet & \bullet & & - \\
\\
\\
\\
\bullet & \bullet \\
& \bullet \bullet & \bullet \bullet & & 2p \\
& \bullet & \bullet \\
\ar@{..}[3,3] & \bullet & \bullet \\
& \bullet & \bullet \\
& & \bullet \bullet \\
& & & \bullet & 0 \\
\\
-3 & -2 & -1 & 0 & s \backslash t
}
$$
\caption{$\bF_p[v_2]$-generators of $E^2_{s,t} = H_{\mot}^{-s}(L_p;
	\bF_{p^2}(t/2)) \Longrightarrow V(1)_{s+t} K(L_p)$}
\label{fig:E2KLp}
\end{figure}

In this case, the motivic truncation property for $L_p$ is perfectly
satisfied, in the sense that
$$
H^r_{\mot}(L_p; \bF_{p^2}(m)) \cong
\begin{cases}
H^r_{\et}(L_p; \bF_{p^2}(m)) & \text{for $0 \le r \le m$,} \\
0 & \text{otherwise,}
\end{cases}
$$
where, by definition,
$$
H^r_{\et}(L_p; \bF_{p^2}(*)) = v_2^{-1} H^r_{\mot}(L_p; \bF_{p^2}(*))
$$
is free over the Laurent polynomial ring $\bF_p[v_2^{\pm1}]$ on the same
generators as in Figure~\ref{fig:E2KLp}.  The (additive)
Euler characteristic
$$
\chi(L_p; \bF_{p^2}(m)) = \sum_{r=0}^3 (-1)^r \dim_{\bF_p}
	H_{\mot}^r(L_p; \bF_{p^2}(m))
$$
is zero for each $m$, cf.~\cite{Tat63}*{Thm.~2.2}.  To the eyes of
algebraic $K$-theory and the hypothetical motivic cohomology, the strictly
commutative ring spectrum~$L_p$ behaves much like a brave new field.
We discuss the role of its (non-commutative) residue ring $L/p$ in the
next section.
\end{example}

The {\'e}tale cohomology of a $p$-adic field~$F$ satisfies local
Tate--Poitou duality \cite{Tat63}*{Thm.~2.1}.  In the case of mod~$p$
coefficients, this is a perfect pairing
$$
H_{\et}^r(F; \bZ/p(m)) \otimes H_{\et}^{2-r}(F; \bZ/p(1-m))
	\overset{\cup}\longto H_{\et}^2(F; \bZ/p(1)) \cong \bZ/p
$$
for each $r$ and~$m$.  For general $p$-power torsion coefficients there
is a perfect pairing taking values in the larger group $H_{\et}^2(F;
\bZ/p^\infty(1)) \cong \bZ/p^\infty$, cf.~\cite{Ser67}*{p.~130}.
The multiplicative structure on $V(1)_* K(L_p)$ is compatible with
an algebra structure on $H_{\mot}^*(L_p; \bF_{p^2}(*))$ such that the
resulting multiplicative structure on $H_{\et}^*(L_p; \bF_{p^2}(*))$
also satisfies \emph{arithmetic duality}.  This can be seen as a
rotational symmetry about $(s, t) = (-3/2, p+1)$ in the variant of
Figure~\ref{fig:E2KLp} where $v_2$ has been inverted.

\begin{conjecture}
For finite extensions $B$ of $L_p$ there is a perfect pairing
$$
H_{\et}^r(B; \bF_{p^2}(m)) \otimes H_{\et}^{3-r}(B; \bF_{p^2}(p+1-m))
        \overset{\cup}\longto H_{\et}^3(B; \bF_{p^2}(p+1)) \cong \bZ/p
$$
for each $r$ and~$m$.
\end{conjecture}

\begin{remark}
The dependence of the twist in $\bF_{p^2}(p+1)$ on the prime~$p$ may
be an artifact of the passage to mod~$p$ and~$v_1$ coefficients.
Let $E_2/(p^\infty, u_1^\infty)$ be the $E_2$-module spectrum defined
by the homotopy cofiber sequences
$E_2 \to p^{-1} E_2 \to E_2/p^\infty$
and
$E_2/p^\infty \to u_1^{-1} E_2/p^\infty \to
	E_2/(p^\infty, u_1^\infty)$.
Its homotopy groups $\pi_* E_2/(p^\infty, u_1^\infty) =
\bW\bF_{p^2}[[u_1]]/(p^\infty, u_1^\infty)[u^{\pm1}]$ are Pontryagin
dual to those of $E_2$.  Then
$$
\bF_{p^2}(p+1) = V(1)_{2p+2} E_2 \cong V(1)_{2p+3}(E_2/p^\infty)
	\cong V(1)_{2p+4} E_2/(p^\infty, u_1^\infty)
$$
and
\begin{align*}
V(1)_{2p+4} E_2/(p^\infty, u_1^\infty)
	&\subset (S/p)_5 E_2/(p^\infty, u_1^\infty) \\
	&\subset \pi_4 E_2/(p^\infty, u_1^\infty)
	= \bW\bF_{p^2}[[u_1]]/(p^\infty, u_1^\infty)(2) \,.
\end{align*}
The conjectured arithmetic duality for mod~$p$ and~$v_1$ coefficients
may be a special case of a duality for $p$- and $u_1$-power torsion
coefficients, taking values in
$$
H^3_{\et}(B; \bW\bF_{p^2}[[u_1]]/(p^\infty, u_1^\infty)(2)) \,.
$$
It would be desirable to find a canonical identification of this group,
like the Hasse invariant in the classical case of $p$-adic fields and
Kato's work \cite{Kat82} on the Galois cohomology of higher-dimensional
local fields.
\end{remark}

\section{Fraction fields and ramified extensions}

The {\'e}tale cohomology of a field is, by construction, the same as its
Galois cohomology, i.e., the continuous group cohomology of its absolute
Galois group.  There is no such direct description of $H_{\et}^r(L_p;
\bF_{p^2}(m))$, since according to Baker--Richter \cite{BR08} the maximal
connected pro-Galois extension of $L_p$ is the composite
$$
L_p \overset{\phi}\longto KU_p \longto KU^{\nr}_p \,,
$$
where $\pi_*(KU^{\nr}_p) = \bW\bar\bF_p[u^{\pm1}]$.  The unramified
extensions of $\pi_0(KU_p) = \bZ_p$ are spectrally realized, using
the methods of Schw{\"a}nzl--Vogt--Waldhausen \cite{SVW99}*{Thm.~3}
or Goerss--Hopkins--Miller \cite{GH04}, but the associated Galois
group only has $p$-cohomo\-logical dimension~$1$, whereas $L_p$ would
have $p$-cohomological dimension~$3$.  Likewise, the maximal connected
pro-Galois extension of $E_n$ is $E_n^{\nr}$, with $\pi_*(E_n^{\nr})
= \bW\bar\bF[[u_1, \dots, u_{n-1}]][u^{\pm1}]$, of $p$-cohomological
dimension~$1$ over~$E_n$ and~$\hat L_n E(n)$.

To allow for ramification at~$p$, one might simply invert that prime.
However, the resulting strictly commutative ring spectrum $p^{-1} L_p$,
with $\pi_*(p^{-1} L_p) = \bQ_p[v_1^{\pm1}]$, is an algebra over $p^{-1}
S_p = H\bQ_p$, so $V(1)_* K(p^{-1} L_p)$ is an algebra over $V(1)_*
K(\bQ_p)$, where $v_2$ acts trivially.  Hence $H^r_{\et}(p^{-1} L_p;
\bF_{p^2}(*))$ would be zero.

A milder form of localization may be appropriate.  By Waldhausen's
localization theorem \cite{Wal85}, the homotopy fiber of $K(L_p) \to
K(p^{-1} L_p)$ is given by the algebraic $K$-theory of the category
with cofibrations of finite cell $L_p$-modules with $p$-power torsion
homotopy, equipped with the usual weak equivalences.  We might instead
step back to a category with cofibrations of coherent $L/p^\nu$-modules
(i.e., having degreewise finite homotopy groups, see Barwick--Lawson
\cite{BL:14}), for some natural number $\nu$, and suppose that these
have the same algebraic $K$-theory as the category with cofibrations
of finite cell $L/p$-modules.  Here $L/p = K(1)$ is the first Morava
$K$-theory, which by Angeltveit \cite{Ang11} is a strict ring spectrum,
but not strictly commutative.  By Davis--Lawson \cite{DL:13}*{Cor.~6.4}
the tower $\{L/p^\nu\}_{\nu}$ is as commutative as possible in the
category of pro-spectra:
$$
L/p \longleftarrow \{L/p^\nu\}_\nu \longleftarrow L_p \longto
	p^{-1} L_p \,.
$$

\begin{definition}
Let $K(\ff L_p)$ be defined by the homotopy cofiber sequence
$$
K(L/p) \overset{i_*}\longto K(L_p) \longto K(\ff L_p) \,,
$$
where $i_*$ is the transfer map associated to $i \colon L_p \to L/p$.
\end{definition}

We think of $K(\ff L_p)$ as the algebraic $K$-theory of a hypothetical
\emph{fraction field} of~$L_p$, intermediate between $L_p$ and $p^{-1}
L_p$, and similar to the $2$-dimensional local field~$\bQ_p(\!(u)\!)$.
Its mod~$p$ and~$v_1$ homotopy groups can be calculated using the
following result, in combination with the homotopy cofiber sequence
$K(\bF_p) \to K(\ell/p) \to K(L/p)$.

\begin{theorem}[Ausoni--Rognes \cite{AR12}*{Thm.~7.6, Thm.~7.7}]
$V(1)_* TC(\ell/p; p)$ is a finitely
generated free $\bF_p[v_2]$-module on $2p^2 - 2p + 8$ generators,
which are located in degrees $-1 \le * \le 2p^2 + 2p - 2$.
There is an exact sequence of $\bF_p[v_2]$-modules
$$
0 \to V(1)_* K(\ell/p) \overset{\trc}\longto
	V(1)_* TC(\ell/p; p) \longto
	\Sigma^{-1} \bF_p \oplus \Sigma^{2p-2} \bF_p \to 0 \,.
$$
\end{theorem}

\begin{example}
The expected $E^2$-term for a spectral sequence
$$
E^2_{s,t} = H_{\mot}^{-s}(\ff L_p; \bF_{p^2}(t/2))
	\Longrightarrow V(1)_{s+t} K(\ff L_p)
$$
is displayed for $p=5$ in Figure~\ref{fig:E2KffLp}.  In addition there
are four sporadic copies of $\bF_p$, in degrees $2p-3$, $2p-2$, $2p-2$
and~$2p-1$.  The motivic truncation properties for $\ff L_p$, analogous
to~\eqref{eq:BLmotet} and \eqref{eq:BLetmot}, are clearly visible,
and conjecturally there is now a perfect arithmetic duality pairing
$$
H_{\et}^r(\ff L_p; \bF_{p^2}(m))
	\otimes H_{\et}^{3-r}(\ff L_p; \bF_{p^2}(2-m))
	\overset{\cup}\longto
	H_{\et}^{3}(\ff L_p; \bF_{p^2}(2)) \cong \bZ/p \,.
$$

\begin{figure}[t]
$$
\xymatrix@!C@!R@R-30pt{
\ar@{..}[3,3] \ar@{..}[24,0] \\
\bullet \\
& \bullet{\bullet}\bullet & & & 2p^2 \\
& \bullet & \bullet \bullet & \ar@{..}[24,0] \\
& \bullet & \bullet \\
& \bullet & \bullet \\
& \bullet & \bullet \\
& \bullet & \bullet & & - \\
& \bullet & \bullet \\
& \bullet & \bullet \\
& \bullet & \bullet \\
& \bullet & \bullet \\
& \bullet & \bullet & & - \\
& \bullet & \bullet \\
& \bullet & \bullet \\
& \bullet & \bullet \\
& \bullet & \bullet \\
& \bullet & \bullet & & - \\
& \bullet & \bullet \\
& \bullet & \bullet \\
& \bullet & \bullet \\
& \bullet & \bullet \\
& \bullet & \bullet & & 2p \\
& \bullet & \bullet \\
\ar@{..}[3,3] & \bullet & \bullet \\
& \bullet \bullet & \bullet \\
& & \bullet{\bullet}\bullet \\
& & & \bullet & 0 \\
\\
-3 & -2 & -1 & 0 & s \backslash t
}
$$
\caption{$\bF_p[v_2]$-generators of $E^2_{s,t} = H_{\mot}^{-s}(\ff L_p;
	\bF_{p^2}(t/2)) \Longrightarrow V(1)_{s+t} K(\ff L_p)$}
\label{fig:E2KffLp}
\end{figure}

After such localization of $L_p$ away from $L/p$, it may be possible
to construct enough Galois extensions of $\ff L_p$ to realize its
$v_2$-localized motivic cohomology as continuous group cohomology
$$
H_{\et}^r(\ff L_p; \bF_{p^2}(m))
	\cong H_{\gp}^r(G_{\ff L_p}; \bF_{p^2}(m)) \,,
$$
for an absolute Galois group $G_{\ff L_p}$ of $p$-cohomological
dimension~$3$, corresponding to some maximal extension $\ff L_p \to
\Omega_1$.  If each Galois extension of $\bQ_p$ can be lifted to an
extension of $\ff L_p$, we get a short exact sequence
$$
1 \to I_{v_1} \longto G_{\ff L_p} \longto G_{\bQ_p} \to 1 \,,
$$
with $I_{v_1}$ the inertia group over~$(v_1)$.  Here $G_{\bQ_p}$
has $p$-cohomological dimension~$2$, and $I_{v_1}$ will have
$p$-cohomological dimension~$1$.
\end{example}

In the less structured setting of ring spectra up to homotopy it is
possible to construct totally ramified extensions of $KU_p$, complementary
to the unramified extension $KU_p^{\nr}$, without inverting~$p$.  Torii
\cite{Tor98}*{Thm.~2.5} shows that for each $r\ge1$ the homotopy cofiber
$$
F(B\bZ/p^{r-1}_+, KU_p) \overset{\tau_r^\#}\longto
	F(B\bZ/p^r_+, KU_p) \longto KU_p[\zeta_{p^r}] \,,
$$
of the map of function spectra induced by the stable transfer map $\tau_r
\colon B\bZ/p^r_+ \to B\bZ/p^{r-1}_+$, is a ring spectrum up to homotopy
with $\pi_* KU_p[\zeta_{p^r}] \cong \bZ_p[\zeta_{p^r}][u^{\pm1}]$,
where $\zeta_{p^r}$ denotes a primitive $p^r$-th root of unity.  (He has
similar realization results in the $K(n)$-local category.)  However,
it does not make sense to talk about the algebraic $K$-theory of a ring
spectrum up to homotopy, so these constructions are only helpful if they
can be made strict.

It is not possible to realize $KU_p[\zeta_{p^r}]$ as a strictly
commutative ring spectrum.  Angeltveit \cite{Ang08}*{Rem.~5.18} uses the
identity $\psi^p(x) = x^p + p \theta(x)$ among power operations in $\pi_0$
of a $K(1)$-local strictly commutative ring spectrum to show that if $-p$
admits a $k$-th root in such a ring, with $k\ge2$, then $p$ is invertible
in that ring.  For $r=1$ we have $\bZ_p[\zeta_p] = \bZ_p[\xi]$ where
$\xi^{p-1} = -p$, so for $p$ odd this proves that adjoining $\zeta_p$
to $KU_p$ in a strictly commutative context will also invert~$p$.

It is, however, possible to adjoin $\zeta_p$ to $\pi_0$ of the connective
cover $ku_p$, in the category of strictly commutative ring spectra,
without fully inverting~$p$.  Instead, one must make some positive power
of the Bott element $u \in \pi_2(ku_p)$ singly divisible by~$p$.  If one
thereafter inverts~$u$, it follows that~$p$ has also become invertible.
To achieve this, we modify Torii's construction for $r=1$ by replacing the
transfer map with a norm map.  This leads to the $G$-Tate construction
$$
B^{tG} = t_G(B)^G = [\widetilde{EG} \wedge F(EG_+, i_*B)]^G
$$
for a spectrum $B$ with $G$-action, cf.~Greenlees--May \cite{GM95}*{p.~3}.
This construction preserves strictly commutative ring structures, see
McClure \cite{McC96}*{Thm.~1}.

\begin{example}
Let
$$
KU'_p[\xi] = (ku_p)^{t\bZ/p}
$$
denote the $\bZ/p$-Tate construction on the spectrum $ku_p$ with trivial
$\bZ/p$-action, and let $ku'_p[\xi] = KU'_p[\xi][0,\infty)$ be its
connective cover.  Additively, these are generalized Eilenberg--Mac\,Lane
spectra, cf.~Davis--Mahowald \cite{DM84} and \cite{GM95}*{Thm.~13.5}.
Multiplicatively, $\pi_*(KU'_p[\xi]) \cong \bZ_p[\xi][v^{\pm1}]$ where
$p + \xi^{p-1} = 0$ and $|v| = 2$.  Furthermore,
$$
\pi_*(ku'_p[\xi]) \cong \bZ_p[\xi][v] \,,
$$
and a morphism $ku_p \to ku'_p[\xi]$ of strictly commutative ring spectra
induces the ring homomorphism $\bZ_p[u] \to \bZ_p[\xi][v]$ that maps
$u$ to $\xi \cdot v$.  There is no multiplicative morphism $KU_p \to
KU'_p[\xi]$, but
$$
KU_p \wedge_{ku_p} ku'_p[\xi] = u^{-1} ku'_p[\xi] \simeq
	KU\bQ_p(\xi) = KU\bQ_p(\zeta_p)
$$
is a totally ramified extension of $KU\bQ_p = p^{-1} KU_p$.  We get
a diagram of strictly commutative ring spectra
$$
\xymatrix{
(H\bZ_p)^{t\bZ/p} & KU'_p[\xi] \ar[l] \ar[r] & KU\bQ_p(\xi) \\
(H\bZ_p)^{t\bZ/p}[0,\infty) \ar[d] \ar[u]
	& ku'_p[\xi] \ar[l] \ar[r] \ar[d] \ar[u]
	& ku\bQ_p(\xi) \ar[d] \ar[u] \\
H\bZ/p & H\bZ_p[\xi] \ar[l] \ar[r] & H\bQ_p(\xi)
}
$$
with horizontal maps reducing modulo or inverting~$\xi$, and vertical maps
reducing modulo or inverting~$v$.  Here $\pi_*((H\bZ_p)^{t\bZ/p}) \cong
\hat H^{-*}(\bZ/p; \bZ_p) = \bZ/p[v^{\pm1}]$.  We view $ku'_p[\xi]$ as an
integral model for a $2$-dimensional local field close to $KU\bQ_p(\xi)$,
but note that $ku'_p[\xi]$ is not finite as a $ku_p$-module.
\end{example}

\begin{example}
Let $(\bZ/p)^\times$ act on the group $\bZ/p$ by multiplication, hence
also on the $\bZ/p$-Tate construction $KU'_p[\xi] = (ku_p)^{t\bZ/p}$.  Let
$$
KU'_p = (KU'_p[\xi])^{h(\bZ/p)^\times}
$$
be the homotopy fixed points, and let $ku'_p = KU'_p[0, \infty)$ be
its connective cover.  These are strictly commutative ring spectra,
with $\pi_*(KU'_p) \cong \bZ_p[u, w^{\pm1}]/(pw+u^{p-1})$ where $|w| =
2p-2$, and
$\pi_*(ku'_p) \cong \bZ_p[u, w]/(pw+u^{p-1})$.
Multiplicative morphisms $ku_p \to ku'_p \to ku'_p[\xi]$ induce the
inclusions $\bZ_p[u] \to \bZ_p[u, w]/(pw+u^{p-1}) \to \bZ_p[\xi][v]$
where $w$ maps to $v^{p-1}$.

The morphism $KU'_p \to KU'_p[\xi]$ is a $(\bZ/p)^\times$-Galois extension
in the sense of \cite{Rog08}.  The morphism $ku'_p \to ku'_p[\xi]$ becomes
$(\bZ/p)^\times$-Galois after inverting~$p$ or~$w$.  It remains to be
determined whether $V(1)_* K(ku'_p)$ remains purely $v_2$-periodic, i.e.,
whether the multiplicative approximation $ku_p \to ku'_p$ counters the
additive splitting of $ku'_p$ as a sum of suspended Eilenberg--Mac\,Lane
spectra.
\end{example}

\begin{example}
More generally, for $r\ge1$ let $G = \bZ/p^r$, let $\sP$ be the
family of proper subgroups of $G$, and let $KU_G[0,\infty)$ be
the ``brutal'' truncation of $G$-equivariant periodic $K$-theory,
cf.~\cite{Gre04}*{p.~129}.  Define
$$
KU'_p[\zeta_{p^r}] = (KU_G[0,\infty))^{t\sP}
	= [\widetilde{E\sP} \wedge F(E\sP_+, KU_G[0,\infty))]^G
$$
to be the $\sP$-Tate construction, as in Greenlees--May
\cite{GM95}*{\S17}, and let $ku'_p[\zeta_{p^r}]$ be its
connective cover.  Then
$\pi_*(KU'_p[\zeta_{p^r}]) \cong \bZ_p[\zeta_{p^r}] [v^{\pm1}]$
and
$$
\pi_*(ku'_p[\zeta_{p^r}]) \cong \bZ_p[\zeta_{p^r}] [v] \,.
$$
The map $(KU_G[0,\infty))^{t\sP} \to (KU_G)^{t\sP}$ induces the
inclusion $\bZ_p[\zeta_{p^r}] \subset \bQ_p[\zeta_{p^r}]$
in each even degree.  To prove this, one can compute Amitsur--Dress
homology for the family $\sP$, use the generalized Tate spectral
sequence from \cite{GM95}*{\S22}, and then compare with the calculation
in \cite{GM95}*{\S19} of the periodic case.  The cyclotomic extension
$\bZ_p[\zeta_{p^r}]$ arises as $(R(G)/J'\!\sP)^\wedge_{J\sP}$, where $R(G)$
is the representation ring, $J\sP$ is the kernel of the restriction map
$R(G) \to R(H)$, and $J'\!\sP$ is the image of the induction map $R(H)
\to R(G)$, where $H$ is the index~$p$ subgroup in $G$.
\end{example}

%

\section{Logarithmic ring spectra}

The heuristics from the last two sections suggest that we should
attempt to construct ramified finite extensions $B \to C$ of strictly
commutative ring spectra $B$ like $\ell_p$, $ku_p$ and $e_n$.  The
Goerss--Hopkins--Miller obstruction theory \cite{GH04} for strictly
commutative $B$-algebra structures on such spectra $C$ has vanishing
obstruction groups in the case of unramified extensions, but appears to be
less useful in the case of ramification over~$(p)$, due to the presence of
nontrivial (topological) Andr{\'e}--Quillen cohomology groups \cite{BR04}.

The same heuristics also suggest that we should approach the
extension problem by passing to mildly local versions of~$B$,
intermediate between $B$ and $p^{-1} B$.  In arithmetic algebraic
geometry, one such intermediary is provided by logarithmic geometry,
cf.~Kato \cite{Kat89}.  An affine pre-log scheme $(\Spec R, M)$ is a
scheme $\Spec R$, a commutative monoid $M$, and a homomorphism $\alpha
\colon M \to (R, \cdot)$ to the underlying multiplicative monoid of~$R$.
More precisely, $M$ and $\alpha$ live {\'e}tale locally on $\Spec R$.
In this wider context, there is a factorization
$$
\Spec R[M^{-1}] \longto (\Spec R, M) \longto \Spec R
$$
of the natural inclusion, and the right-hand map is often a well-behaved
proper replacement for the composite open immersion.  Logarithmic
structures on valuation rings in $p$-adic fields were successfully
used by Hesselholt--Madsen \cite{HM03} to analyze the topological cyclic
homology and algebraic $K$-theory of these classical rings.

A theory of logarithmic structures on strictly commutative ring spectra
was started by the author in \cite{Rog09}, and developed further in
joint work with Sagave and Schlichtkrull.  To present it, we take
the category $\CSp$ of commutative symmetric ring spectra \cite{HSS00},
with the positive stable model structure \cite{MMSS01}*{\S14}, as our
model for strictly commutative ring spectra.

By the graded underlying space of a symmetric spectrum $A$ we mean
a diagram
$$
\Omega^{\cJ}(A) \colon (\bn_1, \bn_2) \longmapsto \Omega^{n_2} A_{n_1}
$$
of spaces, where $(\bn_1, \bn_2)$ ranges over the objects in a
category~$\cJ$.  We call such a diagram a $\cJ$-space.  Following Sagave,
the natural category $\cJ$ to consider turns out to be isomorphic
to Quillen's construction $\Sigma^{-1}\Sigma$, where $\Sigma$ is the
permutative groupoid of finite sets and bijections.  Its nerve $B\cJ \cong
B(\Sigma^{-1} \Sigma)$ is homotopy equivalent to $QS^0 = \Omega^\infty S$.
For any $\cJ$-space~$X$, the homotopy colimit $X_{h\cJ} = \hocolim_{\cJ}
X$ is augmented over $B\cJ$, so we say that $X$ is $QS^0$-graded.  A map
$X \to Y$ of $\cJ$-spaces is called a $\cJ$-equivalence if the induced
map $X_{h\cJ} \to Y_{h\cJ}$ is a weak equivalence.

If $A$ is a commutative symmetric ring spectrum, then $\Omega^{\cJ}(A)$
is a commutative monoid with respect to a convolution product in the
category of $\cJ$-spaces.  The category $\CSJ$ of commutative $\cJ$-space
monoids has a positive projective model structure \cite{SS12}*{\S4},
with the $\cJ$-equivalences as the weak equivalences, and is Quillen
equivalent to a category of $E_\infty$ spaces over $B\cJ$.  The functor
$\Omega^{\cJ} \colon \CSp \to \CSJ$ admits a left adjoint $M \mapsto \SJ[M]$,
and $(\SJ[-], \Omega^{\cJ})$ is a Quillen adjunction.

There is a commutative submonoid of graded homotopy units $\iota \colon
GL_1^{\cJ}(A) \subset \Omega^{\cJ}(A)$.  A \emph{pre-log ring spectrum}
$(A, M, \alpha)$ is a commutative symmetric ring spectrum~$A$ with
a pre-log structure $(M, \alpha)$, i.e., a commutative $\cJ$-space
monoid $M$ and a map $\alpha \colon M \to \Omega^{\cJ}(A)$ in $\CSJ$.
If the pullback $\alpha^{-1}(GL_1^{\cJ}(A)) \to GL_1^{\cJ}(A)$ is a
$\cJ$-equivalence we call $(M, \alpha)$ a log structure and $(A, M,
\alpha)$ a \emph{log ring spectrum}.  We often omit $\alpha$ from the
notation.

In order to classify extensions $(A, M) \to (B, N)$ of pre-log ring
spectra, one is led to study infinitesimal deformations and derivations
in this category.  Derivations are corepresented by a logarithmic version
$T\!AQ(A, M)$ of topological Andr{\'e}--Quillen homology, defined by
a pushout
$$
\xymatrix{
A \wedge_{\SJ[M]} T\!AQ(\SJ[M]) \ar[r]^-{\psi} \ar[d]_-{\alpha}
	& A \wedge \gamma(M) \ar[d] \\
T\!AQ(A) \ar[r]
	& T\!AQ(A, M)
}
$$
of $A$-module spectra, cf.~\cite{Rog09}*{Def.~11.19} and
\cite{Sag14}*{Def.~5.20}.  Here $T\!AQ(A)$ is the ordinary topological
Andr{\'e}--Quillen homology, as defined by Basterra \cite{Bas99}, and
$\gamma(M)$ is the connective spectrum associated to the $E_\infty$ space
$M_{h\cJ}$.  A morphism $(A, M) \to (B, N)$ is formally log {\'e}tale
if $B \wedge_A T\!AQ(A, M) \to T\!AQ(B, N)$ is an equivalence.

Let $j \colon e \to E$ be a fibration of commutative symmetric ring spectra.
The direct image of the trivial log structure on $E$ is the log structure
$(j_* GL_1^{\cJ}(E), \alpha)$ on $e$ given by the pullback
$$
\xymatrix{
j_* GL_1^{\cJ}(E) \ar[r] \ar[d]_-{\alpha}
	& GL_1^{\cJ}(E) \ar[d]^-{\iota} \\
\Omega^{\cJ}(e) \ar[r]^-{\Omega^{\cJ}(j)} & \Omega^{\cJ}(E)
}
$$
in $\CSJ$.
Applying this natural construction to the vertical maps
in~\eqref{eq:SphiHZp}, we get the following example.  Note that $\ell_p
\to ku_p$ is not {\'e}tale, while $L_p \to KU_p$ is $\Delta$-Galois,
hence {\'e}tale.

\begin{theorem}[Sagave \cite{Sag14}*{Thm.~6.1}]
The morphism
$$
\phi \colon (\ell_p, j_* GL_1^{\cJ}(L_p))
	\longto (ku_p, j_* GL_1^{\cJ}(KU_p))
$$
is log {\'e}tale.
\end{theorem}

In order to approximate algebraic $K$-theory, one is likewise led to study
logarithmic topological Hochschild homology and logarithmic topological
cyclic homology.
The former is defined by a pushout
$$
\xymatrix{
\SJ[B^{\cy}(M)] \ar[r]^{\rho} \ar[d]_{\alpha}
	& \SJ[B^{\rep}(M)] \ar[d] \\
THH(A) \ar[r] & THH(A, M)
}
$$
in $\CSp$, cf.~\cite{RSS:I}*{\S4}.  Here $THH(A)$ is the ordinary
topological Hochschild homology of $A$, given by the cyclic
bar construction of $A$ in $\CSp$.  The cyclic bar construction
$B^{\cy}(M)$ is formed in $\CSJ$, and is naturally augmented over $M$.
The \emph{replete bar construction} $B^{\rep}(M)$ can be viewed as
a fibrant replacement of $B^{\cy}(M)$ over~$M$ in a group completion
model structure on $\CSJ$, cf.~\cite{Sag:13}*{Thm.~1.6}, but also has
a more direct description as the (homotopy) pullback in the right hand
square below:
$$
\xymatrix{
B^{\cy}(M) \ar[r]^-{\rho} \ar[d]_{\epsilon}
	& B^{\rep}(M) \ar[r] \ar[d] & B^{\cy}(M^{\gp}) \ar[d]^{\epsilon} \\
M \ar[r]^{=} & M \ar[r]^-{\eta} & M^{\gp} \,.
}
$$
Here $\eta \colon M \to M^{\gp}$ is a group completion in $\CSJ$, which
means that $(M^{\gp})_{h\cJ}$ is a group completion of the $E_\infty$
space $M_{h\cJ}$.  The role of repletion in homotopy theory is similar to
that of working within the subcategory of fine and saturated logarithmic
structures in the discrete setting \cite{Kat89}*{\S2}.

A morphism $(A, M) \to (B, N)$ is formally log thh-{\'e}tale if $B
\wedge_A THH(A, M) \to THH(B, N)$ is an equivalence.  The following
theorem strengthens the previous result.

\begin{theorem}[Rognes--Sagave--Schlichtkrull \cite{RSS:II}*{Thm.~1.5}]
The morphism
$$
\phi \colon (\ell_p, j_* GL_1^{\cJ}(L_p))
	\longto (ku_p, j_* GL_1^{\cJ}(KU_p))
$$
is log thh-{\'e}tale.
\end{theorem}

\begin{remark}
These results harmonize with the classical correspondence between tamely
ramified extensions and log {\'e}tale extensions.  By Noether's theorem
\cite{Noe32}, tame ramification corresponds locally to the existence
of a normal basis.  This conforms with the observation that $ku_p$ is
a retract of a finite cell $\ell_p[\Delta]$-module, so that $\ell_p
\to ku_p$ is tamely ramified.  By contrast, $ku_2$ is not a retract
of a finite cell $ko_2[C_2]$-module, e.g.~because $(ku_2)^{tC_2}$ is
nontrivial, so $ko_2 \to ku_2$ is wildly ramified.
\end{remark}

We say that a commutative symmetric ring spectrum $E$ is $d$-periodic
if $d$ is the minimal positive integer such that $\pi_*(E)$ contains a
unit in degree~$d$.

\begin{theorem}[Rognes--Sagave--Schlichtkrull \cite{RSS:I}*{Thm.~1.5}]
Let $E$ in $\CSp$ be $d$-periodic, with connective cover $j \colon e \to E$.
There is a natural homotopy cofiber sequence
$$
THH(e) \overset{\rho}\longto THH(e, j_* GL_1^{\cJ}(E))
	\overset{\partial}\longto \Sigma THH(e[0,d)) \,,
$$
where $e[0,d)$ is the $(d-1)$-th Postnikov section of $e$.
\end{theorem}

These results allow us to realize the strategy outlined in
\cite{Aus05}*{\S10} to compute the $V(1)$-homotopy of $THH(ku_p)$ by way
of $THH(\ell_p)$, $THH(\ell_p, j_* GL_1^{\cJ}(L_p))$ and $THH(ku_p, j_*
GL_1^{\cJ}(KU_p))$.  The details are given in \cite{RSS:II}*{\S7, \S8}.

When $e[0,d) = H\pi_0(e)$ with $\pi_0(e)$ regular,
Blumberg--Mandell \cite{BlM:14}*{Thm.~4.2.1} have
constructed a map of horizontal homotopy cofiber sequences
$$
\xymatrix{
K(\pi_0(e)) \ar[r]^-{i_*} \ar[d] & K(e) \ar[r]^-{j^*} \ar[d] & K(E) \ar[d] \\
THH(\pi_0(e)) \ar[r]^-{i_*} & THH(e) \ar[r]^-{j^*} & WTHH^\Gamma(e|E)
}
$$
where the vertical arrows are trace maps.

\begin{conjecture}
There is an equivalence of cyclotomic spectra
$$
THH(e, j_* GL_1^{\cJ}(E)) \simeq WTHH^\Gamma(e|E) \,,
$$
compatible with the maps from $THH(e)$ and to $\Sigma THH(\pi_0(e))$.
\end{conjecture}

The author hopes that a logarithmic analog of the Goerss--Hopkins--Miller
obstruction theory \cite{GH04} can be developed to classify log extensions
of log ring spectra, and that in the case of log {\'e}tale extensions the
obstruction groups will vanish in such a way as to enable the construction
of interesting examples.  The underlying strictly commutative ring spectra
should then provide novel examples of tamely ramified extensions, and
realize a larger part of motivic cohomology as a case of Galois cohomology.


\begin{bibdiv}
\begin{biblist}

\bib{Ang08}{article}{
   author={Angeltveit, Vigleik},
   title={Topological Hochschild homology and cohomology of $A_\infty$
   ring spectra},
   journal={Geom. Topol.},
   volume={12},
   date={2008},
   number={2},
   pages={987--1032},
}

\bib{Ang11}{article}{
   author={Angeltveit, Vigleik},
   title={Uniqueness of Morava $K$-theory},
   journal={Compos. Math.},
   volume={147},
   date={2011},
   number={2},
   pages={633--648},
}

\bib{Aus05}{article}{
   author={Ausoni, Christian},
   title={Topological Hochschild homology of connective complex $K$-theory},
   journal={Amer. J. Math.},
   volume={127},
   date={2005},
   number={6},
   pages={1261--1313},
}

\bib{Aus10}{article}{
   author={Ausoni, Christian},
   title={On the algebraic $K$-theory of the complex $K$-theory spectrum},
   journal={Invent. Math.},
   volume={180},
   date={2010},
   number={3},
   pages={611--668},
}

\bib{AR02}{article}{
   author={Ausoni, Christian},
   author={Rognes, John},
   title={Algebraic $K$-theory of topological $K$-theory},
   journal={Acta Math.},
   volume={188},
   date={2002},
   number={1},
   pages={1--39},
}

\bib{AR08}{article}{
   author={Ausoni, Christian},
   author={Rognes, John},
   title={The chromatic red-shift in algebraic K-theory},
   journal={Monographie de L'Enseignement Math{\'e}matique},
   volume={40},
   date={2008},
   pages={13--15},
}

\bib{AR12}{article}{
   author={Ausoni, Christian},
   author={Rognes, John},
   title={Algebraic $K$-theory of the first Morava $K$-theory},
   journal={J. Eur. Math. Soc. (JEMS)},
   volume={14},
   date={2012},
   number={4},
   pages={1041--1079},
}

\bib{BR08}{article}{
   author={Baker, Andrew},
   author={Richter, Birgit},
   title={Galois extensions of Lubin--Tate spectra},
   journal={Homology, Homotopy Appl.},
   volume={10},
   date={2008},
   number={3},
   pages={27--43},
}

\bib{BL:14}{article}{
  author={Barwick, Clark},
  author={Lawson, Tyler},
  title={Regularity of structured ring spectra and localization in K-theory},
  date={2014},
  eprint={ arXiv:1402.6038v2 [math.KT] },
}

\bib{Bas99}{article}{
   author={Basterra, M.},
   title={Andr\'e-Quillen cohomology of commutative $S$-algebras},
   journal={J. Pure Appl. Algebra},
   volume={144},
   date={1999},
   number={2},
   pages={111--143},
}

\bib{BR04}{article}{
   author={Basterra, Maria},
   author={Richter, Birgit},
   title={(Co-)homology theories for commutative ($S$-)algebras},
   conference={
      title={Structured ring spectra},
   },
   book={
      series={London Math. Soc. Lecture Note Ser.},
      volume={315},
      publisher={Cambridge Univ. Press},
      place={Cambridge},
   },
   date={2004},
   pages={115--131},
}

\bib{BeM:14}{article}{
   author={Berglund, A.},
   author={Madsen, I.},
   title={Rational homotopy theory of automorphisms of highly connected
	manifolds},
   date={2014},
   eprint={ arXiv:1401.4096 [math.AT] },
}

\bib{Blo86}{article}{
   author={Bloch, Spencer},
   title={Algebraic cycles and higher $K$-theory},
   journal={Adv. in Math.},
   volume={61},
   date={1986},
   number={3},
   pages={267--304},
}

\bib{BM08}{article}{
   author={Blumberg, Andrew J.},
   author={Mandell, Michael A.},
   title={The localization sequence for the algebraic $K$-theory of
   topological $K$-theory},
   journal={Acta Math.},
   volume={200},
   date={2008},
   number={2},
   pages={155--179},
}

\bib{BlM:14}{article}{
   author={Blumberg, Andrew J.},
   author={Mandell, Michael A.},
   title={Localization for $THH(ku)$ and the topological Hochschild and
	cyclic homology of Waldhausen categories},
   date={2014},
   eprint={ arXiv:1111.4003v3 [math.KT] },
}

\bib{BHM93}{article}{
   author={B{\"o}kstedt, M.},
   author={Hsiang, W. C.},
   author={Madsen, I.},
   title={The cyclotomic trace and algebraic $K$-theory of spaces},
   journal={Invent. Math.},
   volume={111},
   date={1993},
   number={3},
   pages={465--539},
}

\bib{BM94}{article}{
   author={B{\"o}kstedt, M.},
   author={Madsen, I.},
   title={Topological cyclic homology of the integers},
   note={$K$-theory (Strasbourg, 1992)},
   journal={Ast\'erisque},
   number={226},
   date={1994},
   pages={7--8, 57--143},
}

\bib{BM95}{article}{
   author={B{\"o}kstedt, M.},
   author={Madsen, I.},
   title={Algebraic $K$-theory of local number fields: the unramified case},
   conference={
      title={Prospects in topology},
      address={Princeton, NJ},
      date={1994},
   },
   book={
      series={Ann. of Math. Stud.},
      volume={138},
      publisher={Princeton Univ. Press},
      place={Princeton, NJ},
   },
   date={1995},
   pages={28--57},
}

\bib{BW87}{article}{
   author={B{\"o}kstedt, Marcel},
   author={Waldhausen, Friedhelm},
   title={The map $BSG\to A(*)\to QS^0$},
   conference={
      title={Algebraic topology and algebraic $K$-theory (Princeton, N.J.,
      1983)},
   },
   book={
      series={Ann. of Math. Stud.},
      volume={113},
      publisher={Princeton Univ. Press},
      place={Princeton, NJ},
   },
   date={1987},
   pages={418--431},
}

\bib{Bor74}{article}{
   author={Borel, Armand},
   title={Stable real cohomology of arithmetic groups},
   journal={Ann. Sci. \'Ecole Norm. Sup. (4)},
   volume={7},
   date={1974},
   pages={235--272 (1975)},
}

\bib{Bro72}{book}{
   author={Browder, William},
   title={Surgery on simply-connected manifolds},
   note={Ergebnisse der Mathematik und ihrer Grenzgebiete, Band 65},
   publisher={Springer-Verlag},
   place={New York},
   date={1972},
   pages={ix+132},
}


\bib{DL:13}{article}{
   author={Davis, Daniel G.},
   author={Lawson, Tyler},
   title={Commutative ring objects in pro-categories and generalized
   	Moore spectra},
   date={2013},
   eprint={ arXiv:1208.4519v3 [math.AT] },
}

\bib{DM84}{article}{
   author={Davis, Donald M.},
   author={Mahowald, Mark},
   title={The spectrum $(P\wedge b{\rm o})_{-\infty }$},
   journal={Math. Proc. Cambridge Philos. Soc.},
   volume={96},
   date={1984},
   number={1},
   pages={85--93},
}

\bib{Dun97}{article}{
   author={Dundas, Bj{\o}rn Ian},
   title={Relative $K$-theory and topological cyclic homology},
   journal={Acta Math.},
   volume={179},
   date={1997},
   number={2},
   pages={223--242},
}

\bib{DM98}{article}{
   author={Dwyer, W. G.},
   author={Mitchell, S. A.},
   title={On the $K$-theory spectrum of a ring of algebraic integers},
   journal={$K$-Theory},
   volume={14},
   date={1998},
   number={3},
   pages={201--263},
}

\bib{EKMM97}{book}{
   author={Elmendorf, A. D.},
   author={Kriz, I.},
   author={Mandell, M. A.},
   author={May, J. P.},
   title={Rings, modules, and algebras in stable homotopy theory},
   series={Mathematical Surveys and Monographs},
   volume={47},
   note={With an appendix by M. Cole},
   publisher={American Mathematical Society},
   place={Providence, RI},
   date={1997},
   pages={xii+249},
}

\bib{FH78}{article}{
   author={Farrell, F. T.},
   author={Hsiang, W. C.},
   title={On the rational homotopy groups of the diffeomorphism groups of
   discs, spheres and aspherical manifolds},
   conference={
      title={Algebraic and geometric topology (Proc. Sympos. Pure Math.,
      Stanford Univ., Stanford, Calif., 1976), Part 1},
   },
   book={
      series={Proc. Sympos. Pure Math., XXXII},
      publisher={Amer. Math. Soc.},
      place={Providence, R.I.},
   },
   date={1978},
   pages={325--337},
}

\bib{FJ91}{article}{
   author={Farrell, F. T.},
   author={Jones, L. E.},
   title={Stable pseudoisotopy spaces of compact non-positively curved
   manifolds},
   journal={J. Differential Geom.},
   volume={34},
   date={1991},
   number={3},
   pages={769--834},
}

\bib{GRW:12}{article}{
  author={Galatius, S.},
  author={Randal-Williams, O.},
  title={Stable moduli spaces of high dimensional manifolds},
  date={2012},
  eprint={ arXiv:1201.3527v2 [math.AT] },
}

\bib{GL01}{article}{
   author={Geisser, Thomas},
   author={Levine, Marc},
   title={The Bloch--Kato conjecture and a theorem of Suslin--Voevodsky},
   journal={J. Reine Angew. Math.},
   volume={530},
   date={2001},
   pages={55--103},
}

\bib{GH04}{article}{
   author={Goerss, P. G.},
   author={Hopkins, M. J.},
   title={Moduli spaces of commutative ring spectra},
   conference={
      title={Structured ring spectra},
   },
   book={
      series={London Math. Soc. Lecture Note Ser.},
      volume={315},
      publisher={Cambridge Univ. Press},
      place={Cambridge},
   },
   date={2004},
   pages={151--200},
}

\bib{Goo86}{article}{
   author={Goodwillie, Thomas G.},
   title={Relative algebraic $K$-theory and cyclic homology},
   journal={Ann. of Math. (2)},
   volume={124},
   date={1986},
   number={2},
   pages={347--402},
}

\bib{Goo91}{article}{
   author={Goodwillie, Thomas G.},
   title={The differential calculus of homotopy functors},
   conference={
      title={Proceedings of the International Congress of
              Mathematicians, Vol.\ II},
      address={Kyoto},
      date={1990},
   },
   book={
      publisher={Math. Soc. Japan},
      place={Tokyo},
   },
   date={1991},
   pages={621--630},
}

\bib{Gre04}{article}{
   author={Greenlees, J. P. C.},
   title={Equivariant connective $K$-theory for compact Lie groups},
   journal={J. Pure Appl. Algebra},
   volume={187},
   date={2004},
   number={1-3},
   pages={129--152},
}

\bib{GM95}{article}{
   author={Greenlees, J. P. C.},
   author={May, J. P.},
   title={Generalized Tate cohomology},
   journal={Mem. Amer. Math. Soc.},
   volume={113},
   date={1995},
   number={543},
   pages={viii+178},
}

\bib{Hat75}{article}{
   author={Hatcher, A. E.},
   title={Higher simple homotopy theory},
   journal={Ann. of Math. (2)},
   volume={102},
   date={1975},
   number={1},
   pages={101--137},
}


\bib{Hes09}{article}{
   author={Hesselholt, Lars},
   title={On the Whitehead spectrum of the circle},
   conference={
      title={Algebraic topology},
   },
   book={
      series={Abel Symp.},
      volume={4},
      publisher={Springer},
      place={Berlin},
   },
   date={2009},
   pages={131--184},
}

\bib{HM03}{article}{
   author={Hesselholt, Lars},
   author={Madsen, Ib},
   title={On the $K$-theory of local fields},
   journal={Ann. of Math. (2)},
   volume={158},
   date={2003},
   number={1},
   pages={1--113},
}

\bib{HS98}{article}{
   author={Hopkins, Michael J.},
   author={Smith, Jeffrey H.},
   title={Nilpotence and stable homotopy theory. II},
   journal={Ann. of Math. (2)},
   volume={148},
   date={1998},
   number={1},
   pages={1--49},
}

\bib{HSS00}{article}{
   author={Hovey, Mark},
   author={Shipley, Brooke},
   author={Smith, Jeff},
   title={Symmetric spectra},
   journal={J. Amer. Math. Soc.},
   volume={13},
   date={2000},
   number={1},
   pages={149--208},
}

\bib{HJ82}{article}{
   author={Hsiang, W. C.},
   author={Jahren, B.},
   title={A note on the homotopy groups of the diffeomorphism groups of
   spherical space forms},
   conference={
      title={Algebraic $K$-theory, Part II},
      address={Oberwolfach},
      date={1980},
   },
   book={
      series={Lecture Notes in Math.},
      volume={967},
      publisher={Springer},
      place={Berlin},
   },
   date={1982},
   pages={132--145},
}

\bib{Igu88}{article}{
   author={Igusa, Kiyoshi},
   title={The stability theorem for smooth pseudoisotopies},
   journal={$K$-Theory},
   volume={2},
   date={1988},
   number={1-2},
   pages={vi+355},
}

\bib{Kat82}{article}{
   author={Kato, Kazuya},
   title={Galois cohomology of complete discrete valuation fields},
   conference={
      title={Algebraic $K$-theory, Part II},
      address={Oberwolfach},
      date={1980},
   },
   book={
      series={Lecture Notes in Math.},
      volume={967},
      publisher={Springer},
      place={Berlin},
   },
   date={1982},
   pages={215--238},
}

\bib{Kat89}{article}{
   author={Kato, Kazuya},
   title={Logarithmic structures of Fontaine--Illusie},
   conference={
      title={Algebraic analysis, geometry, and number theory (Baltimore, MD,
      1988)},
   },
   book={
      publisher={Johns Hopkins Univ. Press},
      place={Baltimore, MD},
   },
   date={1989},
   pages={191--224},
}

\bib{KR97}{article}{
   author={Klein, John R.},
   author={Rognes, John},
   title={The fiber of the linearization map $A(*)\to K({\bf Z})$},
   journal={Topology},
   volume={36},
   date={1997},
   number={4},
   pages={829--848},
}

\bib{MW07}{article}{
   author={Madsen, Ib},
   author={Weiss, Michael},
   title={The stable moduli space of Riemann surfaces: Mumford's conjecture},
   journal={Ann. of Math. (2)},
   volume={165},
   date={2007},
   number={3},
   pages={843--941},
}

\bib{MP89}{collection}{
   title={Algebraic topology},
   series={Contemporary Mathematics},
   volume={96},
   booktitle={Proceedings of the International Conference held at
   Northwestern University, Evanston, Illinois, March 21--24, 1988},
   editor={Mahowald, Mark},
   editor={Priddy, Stewart},
   publisher={American Mathematical Society},
   place={Providence, RI},
   date={1989},
   pages={xiv+350},
}

\bib{MMSS01}{article}{
   author={Mandell, M. A.},
   author={May, J. P.},
   author={Schwede, S.},
   author={Shipley, B.},
   title={Model categories of diagram spectra},
   journal={Proc. London Math. Soc. (3)},
   volume={82},
   date={2001},
   number={2},
   pages={441--512},
}

\bib{McC97}{article}{
   author={McCarthy, Randy},
   title={Relative algebraic $K$-theory and topological cyclic homology},
   journal={Acta Math.},
   volume={179},
   date={1997},
   number={2},
   pages={197--222},
}

\bib{McC96}{article}{
   author={McClure, J. E.},
   title={$E_\infty$-ring structures for Tate spectra},
   journal={Proc. Amer. Math. Soc.},
   volume={124},
   date={1996},
   number={6},
   pages={1917--1922},
}

\bib{MS93}{article}{
   author={McClure, J. E.},
   author={Staffeldt, R. E.},
   title={On the topological Hochschild homology of $b{\rm u}$. I},
   journal={Amer. J. Math.},
   volume={115},
   date={1993},
   number={1},
   pages={1--45},
}

\bib{Mos68}{article}{
   author={Mosher, Robert E.},
   title={Some stable homotopy of complex projective space},
   journal={Topology},
   volume={7},
   date={1968},
   pages={179--193},
}

\bib{Noe32}{article}{
   author={Noether, E.},
   title={Normalbasis bei K\"orpern ohne h\"ohere Verzweigung},
   journal={J. Reine Angew. Math.},
   volume={167},
   pages={147--152},
   date={1932},
   publisher={Walter de Gruyter, Berlin},
}

\bib{Qui73}{article}{
   author={Quillen, Daniel},
   title={Higher algebraic $K$-theory. I},
   conference={
      title={Algebraic $K$-theory, I: Higher $K$-theories (Proc. Conf.,
      Battelle Memorial Inst., Seattle, Wash., 1972)},
   },
   book={
      publisher={Springer},
      place={Berlin},
   },
   date={1973},
}

\bib{Qui75}{article}{
   author={Quillen, Daniel},
   title={Higher algebraic $K$-theory},
   conference={
      title={Proceedings of the International Congress of Mathematicians
      (Vancouver, B. C., 1974), Vol. 1},
   },
   book={
      publisher={Canad. Math. Congress, Montreal, Que.},
   },
   date={1975},
   pages={171--176},
}

\bib{Rog93}{article}{
   author={Rognes, John},
   title={Characterizing connected $K$-theory by homotopy groups},
   journal={Math. Proc. Cambridge Philos. Soc.},
   volume={114},
   date={1993},
   number={1},
   pages={99--102},
}

\bib{Rog98}{article}{
   author={Rognes, John},
   title={Trace maps from the algebraic $K$-theory of the integers (after
   Marcel B\"okstedt)},
   journal={J. Pure Appl. Algebra},
   volume={125},
   date={1998},
   number={1-3},
   pages={277--286},
}

\bib{Rog99a}{article}{
   author={Rognes, John},
   title={The product on topological Hochschild homology of the integers
   with mod $4$ coefficients},
   journal={J. Pure Appl. Algebra},
   volume={134},
   date={1999},
   number={3},
   pages={211--218},
}

\bib{Rog99b}{article}{
   author={Rognes, John},
   title={Topological cyclic homology of the integers at two},
   journal={J. Pure Appl. Algebra},
   volume={134},
   date={1999},
   number={3},
   pages={219--286},
}

\bib{Rog99c}{article}{
   author={Rognes, John},
   title={Algebraic $K$-theory of the two-adic integers},
   journal={J. Pure Appl. Algebra},
   volume={134},
   date={1999},
   number={3},
   pages={287--326},
}

\bib{Rog02}{article}{
   author={Rognes, John},
   title={Two-primary algebraic $K$-theory of pointed spaces},
   journal={Topology},
   volume={41},
   date={2002},
   number={5},
   pages={873--926},
}

\bib{Rog03}{article}{
   author={Rognes, John},
   title={The smooth Whitehead spectrum of a point at odd regular primes},
   journal={Geom. Topol.},
   volume={7},
   date={2003},
   pages={155--184},
}

\bib{Rog08}{article}{
   author={Rognes, John},
   title={Galois extensions of structured ring spectra. Stably dualizable
   groups},
   journal={Mem. Amer. Math. Soc.},
   volume={192},
   date={2008},
   number={898},
   pages={viii+137},
}

\bib{Rog09}{article}{
   author={Rognes, John},
   title={Topological logarithmic structures},
   conference={
      title={New topological contexts for Galois theory and algebraic
      geometry (BIRS 2008)},
   },
   book={
      series={Geom. Topol. Monogr.},
      volume={16},
      publisher={Geom. Topol. Publ., Coventry},
   },
   date={2009},
   pages={401--544},
}

\bib{RSS:I}{article}{
  author={Rognes, John},
  author={Sagave, Steffen},
  author={Schlichtkrull, Christian},
  title={Localization sequences for logarithmic topological Hochschild
  	homology},
  date={2014},
  eprint={ arXiv:1402.1317 [math.AT] },
}

\bib{RSS:II}{article}{
  author={Rognes, John},
  author={Sagave, Steffen},
  author={Schlichtkrull, Christian},
  title={Logarithmic topological Hochschild homology of topological $K$-theory spectra},
  date={2014},
}

\bib{RW00}{article}{
   author={Rognes, J.},
   author={Weibel, C.},
   title={Two-primary algebraic $K$-theory of rings of integers in number
   fields},
   note={Appendix A by Manfred Kolster},
   journal={J. Amer. Math. Soc.},
   volume={13},
   date={2000},
   number={1},
   pages={1--54},
}

\bib{Sag:13}{article}{
  author={Sagave, Steffen},
  title={Spectra of units for periodic ring spectra and group completion
	of graded $E_\infty$ spaces},
  date={2013},
  eprint={ arXiv:1111.6731v2 [math.AT] },
}

\bib{Sag14}{article}{
   author={Sagave, Steffen},
   title={Logarithmic structures on topological $K$-theory spectra},
   journal={Geom. Topol.},
   volume={18},
   date={2014},
   number={1},
   pages={447--490},
}

\bib{SS12}{article}{
   author={Sagave, Steffen},
   author={Schlichtkrull, Christian},
   title={Diagram spaces and symmetric spectra},
   journal={Adv. Math.},
   volume={231},
   date={2012},
   number={3-4},
   pages={2116--2193},
}

\bib{SVW99}{article}{
   author={Schw{\"a}nzl, R.},
   author={Vogt, R. M.},
   author={Waldhausen, F.},
   title={Adjoining roots of unity to $E_\infty$ ring spectra in good
   cases---a remark},
   conference={
      title={Homotopy invariant algebraic structures},
      address={Baltimore, MD},
      date={1998},
   },
   book={
      series={Contemp. Math.},
      volume={239},
      publisher={Amer. Math. Soc.},
      place={Providence, RI},
   },
   date={1999},
   pages={245--249},
}

\bib{Ser67}{article}{
   author={Serre, J.-P.},
   title={Local class field theory},
   conference={
      title={Algebraic Number Theory (Proc. Instructional Conf., Brighton,
      1965)},
   },
   book={
      publisher={Thompson, Washington, D.C.},
   },
   date={1967},
   pages={128--161},
}

\bib{Sma62}{article}{
   author={Smale, S.},
   title={On the structure of manifolds},
   journal={Amer. J. Math.},
   volume={84},
   date={1962},
   pages={387--399},
}

\bib{Sou81}{article}{
   author={Soul{\'e}, Christophe},
   title={On higher $p$-adic regulators},
   conference={
      title={Algebraic $K$-theory, Evanston 1980 (Proc. Conf., Northwestern
      Univ., Evanston, Ill., 1980)},
   },
   book={
      series={Lecture Notes in Math.},
      volume={854},
      publisher={Springer},
      place={Berlin},
   },
   date={1981},
   pages={372--401},
}

\bib{SV00}{article}{
   author={Suslin, Andrei},
   author={Voevodsky, Vladimir},
   title={Bloch--Kato conjecture and motivic cohomology with finite
   coefficients},
   conference={
      title={The arithmetic and geometry of algebraic cycles},
      address={Banff, AB},
      date={1998},
   },
   book={
      series={NATO Sci. Ser. C Math. Phys. Sci.},
      volume={548},
      publisher={Kluwer Acad. Publ.},
      place={Dordrecht},
   },
   date={2000},
   pages={117--189},
}

\bib{Tat63}{article}{
   author={Tate, John},
   title={Duality theorems in Galois cohomology over number fields},
   conference={
      title={Proc. Internat. Congr. Mathematicians},
      address={Stockholm},
      date={1962},
   },
   book={
      publisher={Inst. Mittag-Leffler},
      place={Djursholm},
   },
   date={1963},
   pages={288--295},
}

\bib{Tor98}{article}{
   author={Torii, Takeshi},
   title={Topological realization of the integer ring of local field},
   journal={J. Math. Kyoto Univ.},
   volume={38},
   date={1998},
   number={4},
   pages={781--788},
}

\bib{Voe03}{article}{
   author={Voevodsky, Vladimir},
   title={Motivic cohomology with $\bZ/2$-coefficients},
   journal={Publ. Math. Inst. Hautes \'Etudes Sci.},
   number={98},
   date={2003},
   pages={59--104},
}

\bib{Voe11}{article}{
   author={Voevodsky, Vladimir},
   title={On motivic cohomology with $\bZ/l$-coefficients},
   journal={Ann. of Math. (2)},
   volume={174},
   date={2011},
   number={1},
   pages={401--438},
}

\bib{Wal78}{article}{
   author={Waldhausen, Friedhelm},
   title={Algebraic $K$-theory of topological spaces. I},
   conference={
      title={Algebraic and geometric topology (Proc. Sympos. Pure Math.,
      Stanford Univ., Stanford, Calif., 1976), Part 1},
   },
   book={
      series={Proc. Sympos. Pure Math., XXXII},
      publisher={Amer. Math. Soc.},
      place={Providence, R.I.},
   },
   date={1978},
   pages={35--60},
}

\bib{Wal84}{article}{
   author={Waldhausen, Friedhelm},
   title={Algebraic $K$-theory of spaces, localization, and the chromatic
   filtration of stable homotopy},
   conference={
      title={Algebraic topology, Aarhus 1982},
      address={Aarhus},
      date={1982},
   },
   book={
      series={Lecture Notes in Math.},
      volume={1051},
      publisher={Springer},
      place={Berlin},
   },
   date={1984},
   pages={173--195},
}

\bib{Wal85}{article}{
   author={Waldhausen, Friedhelm},
   title={Algebraic $K$-theory of spaces},
   conference={
      title={Algebraic and geometric topology},
      address={New Brunswick, N.J.},
      date={1983},
   },
   book={
      series={Lecture Notes in Math.},
      volume={1126},
      publisher={Springer},
      place={Berlin},
   },
   date={1985},
   pages={318--419},
}

\bib{WJR13}{book}{
   author={Waldhausen, Friedhelm},
   author={Jahren, Bj{\o}rn},
   author={Rognes, John},
   title={Spaces of PL manifolds and categories of simple maps},
   series={Annals of Mathematics Studies},
   volume={186},
   publisher={Princeton University Press},
   place={Princeton, NJ},
   date={2013},
}

\bib{Wll70}{book}{
   author={Wall, C. T. C.},
   title={Surgery on compact manifolds},
   note={London Mathematical Society Monographs, No. 1},
   publisher={Academic Press},
   place={London},
   date={1970},
   pages={x+280},
}

\bib{WW88}{article}{
   author={Weiss, Michael},
   author={Williams, Bruce},
   title={Automorphisms of manifolds and algebraic $K$-theory. I},
   journal={$K$-Theory},
   volume={1},
   date={1988},
   number={6},
   pages={575--626},
}

\bib{WW01}{article}{
   author={Weiss, Michael},
   author={Williams, Bruce},
   title={Automorphisms of manifolds},
   conference={
      title={Surveys on surgery theory, Vol. 2},
   },
   book={
      series={Ann. of Math. Stud.},
      volume={149},
      publisher={Princeton Univ. Press},
      place={Princeton, NJ},
   },
   date={2001},
   pages={165--220},
}


\end{biblist}
\end{bibdiv}
\end{document}